\documentclass[journal,comsoc]{IEEEtran}
\usepackage[T1]{fontenc}% optional T1 font encoding

\usepackage{graphicx}
\usepackage{color}
\usepackage{cite}
\usepackage{lipsum,graphicx,subcaption}
\usepackage{amsfonts}
\usepackage{mathtools}

\usepackage{algorithm}% http://ctan.org/pkg/algorithm
\usepackage{algpseudocode}% http://ctan.org/pkg/algorithmicx

\ifCLASSINFOpdf
\else
\fi
\usepackage{enumerate}
\usepackage{float}
\usepackage{graphicx}
\usepackage[utf8]{inputenc}
\usepackage{array}
\usepackage{wrapfig}
\usepackage{multirow}
\usepackage{tabu}
\usepackage{algorithm, algpseudocode}
\usepackage{amsmath, amssymb}
\usepackage[utf8]{inputenc}
\usepackage[T1]{fontenc}
\usepackage{textcomp}
\usepackage{siunitx}
\usepackage{color}
\usepackage{multicol}

\usepackage{mathtools}

\usepackage[cmintegrals]{newtxmath}
\begin{document}
%
% paper title
% Titles are generally capitalized except for words such as a, an, and, as,
% at, but, by, for, in, nor, of, on, or, the, to and up, which are usually
% not capitalized unless they are the first or last word of the title.
% Linebreaks \\ can be used within to get better formatting as desired.
% Do not put math or special symbols in the title.
\title{Voltage Feasibility Boundaries for Power System Security Assessment}
%

%Computation of the Voltage Feasibility Boundaries of the Power Flow Solution Space.

% author names and IEEE memberships
% note positions of commas and nonbreaking spaces (  ) LaTeX will not break
% a structure at a  so this keeps an author's name from being broken across
% two lines.
% use \thanks{} to gain access to the first footnote area
% a separate \thanks must be used for each paragraph as LaTeX2e's \thanks
% was not built to handle multiple paragraphs
%

\author{Mazhar~Ali,
        Elena~Gryazina,
        Anatoly~Dymarsky,
        and~Petr~Vorobev,~\IEEEmembership{Member~IEEE}}% <-this % stops a space
%\thanks{Mazhar~Ali, E.Gryazina and Petr~Vorobev are with the Skolkovo Institute of Science \& Technology, Moscow, Russia. A. Dymarsky is affiliated with the Skolkovo Institute of Science \& Technology and the University of Kentucky, Lexington, USA.}}% <-this % stops a space
%\thanks{J. Doe and J. Doe are with Anonymous University.}% <-this % stops a space
%\thanks{Manuscript received April 19, 2005; revised August 26, 2015.}}

% note the % following the last \IEEEmembership and also \thanks - 
% these prevent an unwanted space from occurring between the last author name
% and the end of the author line. i.e., if you had this:
% 
% \author{....lastname \thanks{...} \thanks{...} }
%                     ^------------^------------^----Do not want these spaces!
%
% a space would be appended to the last name and could cause every name on that
% line to be shifted left slightly. This is one of those "LaTeX things". For
% instance, "\textbf{A} \textbf{B}" will typeset as "A B" not "AB". To get
% "AB" then you have to do: "\textbf{A}\textbf{B}"
% \thanks is no different in this regard, so shield the last } of each \thanks
% that ends a line with a % and do not let a space in before the next \thanks.
% Spaces after \IEEEmembership other than the last one are OK (and needed) as
% you are supposed to have spaces between the names. For what it is worth,
% this is a minor point as most people would not even notice if the said evil
% space somehow managed to creep in.

% The paper headers
\markboth{This work has been submitted to IEEE Transactions on Power Systems for possible publication.}%
{Shell \MakeLowercase{\textit{et al.}}: Bare Demo of IEEEtran.cls for IEEE Communications Society Journals}
% The only time the second header will appear is for the odd numbered pages
% after the title page when using the twoside option.
% 
% *** Note that you probably will NOT want to include the author's ***
% *** name in the headers of peer review papers.                   ***
% You can use \ifCLASSOPTIONpeerreview for conditional compilation here if
% you desire.

% If you want to put a publisher's ID mark on the page you can do it like
% this:
%\IEEEpubid{0000--0000/00\$00.00~\copyright~2015 IEEE}
% Remember, if you use this you must call \IEEEpubidadjcol in the second
% column for its text to clear the IEEEpubid mark.

% use for special paper notices
%\IEEEspecialpapernotice{(Invited Paper)}

\renewcommand{\baselinestretch}{.9}

% make the title area
\maketitle

% As a general rule, do not put math, special symbols or citations
% in the abstract or keywords.
\begin{abstract}
Modern power systems face a grand challenge in grid management due to increased electricity demand, imminent disturbances, and uncertainties associated with renewable generation, which can compromise grid security. The security assessment is directly connected to the robustness of the operating condition  and is evaluated by analyzing proximity to the power flow solution space's boundary. Calculating location of such a boundary is a computationally challenging task, linked to the power flow equations' non-linear nature, presence of technological constraints, and complicated network topology. In this paper we introduce a general framework to characterize points on the power flow solution space boundary in terms of auxiliary variables subject to algebraic constraints. Then we develop  an adaptive continuation algorithm to trace 1-dimensional sections of boundary curves which exhibits robust performance and computational tractability. Implementation of the algorithm is described in detail, and its performance is validated on different test networks.
\end{abstract}


% For peer review papers, you can put extra information on the cover
% page as needed:
% \ifCLASSOPTIONpeerreview
% \begin{center} \bfseries EDICS Category: 3-BBND \end{center}
% \fi
%
% For peerreview papers, this IEEEtran command inserts a page break and
% creates the second title. It will be ignored for other modes.
\IEEEpeerreviewmaketitle

\section{Introduction}
Power system security has always been a subject of paramount importance as electric grids transform into convoluted structures to address the growing electricity demand \cite{morison2004power}. Furthermore, integration of Distributed Energy Resources (DERs) like wind turbines, photovoltaics, or other renewables has introduced an abstruse uncertainty in the grid's operation, begetting tremendous concern for its stability. Hence, the security assessment is crucial to ensure secure and economical operation by evaluating operational state robustness against potential contingencies \cite{morison2004power}. It has, consequently, been approached by analyzing the boundaries of the power flow solution map.

Power flow solution boundaries separate regions in the parameter space corresponding to the existence of real-valued solution of the power flow equations. \cite{ajjarapu2007computational}. Quantifying boundary curves is indispensable for planning and operational purposes to reduce the risk of disruption, optimize the power transfers, and take full advantage of the transmission capabilities \cite{mittelstaedt2016iterative}.  Also, to expedite the decision-making process to maintain the power system security and price stability \cite{Hiskens:2001ie,lesieutre2005convexity}.

Mathematically, points on the solution space boundary satisfy the power flow equations' real-valued solution, and the degenerate condition of the power flow Jacobian (referred to as ``Transversality'' condition) \cite{ajjarapu2007computational, Hiskens:2001ie}. These conditions describe the boundary of solvability region. To calculate the feasibility space boundary, the operational inequalities like limits on voltage setpoints, transmission line thermal limits, and generator power outputs, etc., must be valid along with the conditions stated earlier. This paper aims to investigate the feasibility boundary of the power flow equations.

In early studies, the power flow solution space structure was hypothesized as convex \cite{jarjis1981quantitative}. However, an example from \cite{Hiskens:2001ie} represents a non-convex topology with inner folds, manifesting multiple equilibria regions. This embedded complexity within the power flow solution space makes calculating its boundary a challenging task from the computational context. The past contributions fall into two main categories. The first type focuses on calculating only a single point on the power flow solution boundary. They include Continuation Power Flow (CPF) approaches \cite{ajjarapu1992continuation, karystianos2007maximizing}, direct methods from \cite{ ali2017transversality}, nonlinear programming techniques from   \cite{van1991method, ayuev2016fast}, and some recently published innovative algorithms as described in \cite{rao2017estimating}.

The second category of efforts investigates the power flow solution's boundary in the multi-dimensional space. Only a handful of steps are made in this direction. The most notable work from \cite{Hiskens:2001ie} introduces a Predictor-Corrector (PC) continuation method to compute the solution space boundary nomograms. It uses tangents as predictors and hyperplanes as correctors. Such homotopy path tracking routines have a difficult implementation and exhibit slow convergence for tracing curves with sharp edges, or non-convexities  \cite{yamamura1993simple}.  Also, the formulation in \cite{Hiskens:2001ie} uses a right eigenvector transversality condition to enforce the non-trivial kernel of the power flow Jacobian, and it doubles the number of system variables. Hence, resulting in an additional computational burden. This condition also makes the initialization of the continuation procedure sensitive to a good guess of the right eigenvector. Moreover, the formulation from \cite{Hiskens:2001ie} does not consider any operational inequalities. Therefore, it is confined to study only the solvability space of the power flow equations.

An alternative approach to study the boundary of the power flow solution space was made by the Numerical Polynomial Homotopy Continuation (NPHC) method. The contributions from \cite{chandra2017locating} and \cite{molzahn2017computing} show an implementation of NPHC method to compute the power flow equations' feasibility space. The inherent mathematical framework of the NPHC algorithm guarantees to find all isolated solutions of the power flow equations. It uses a PC continuation scheme to trace all these solutions along with a one-dimensional parametrization. Despite the computational robustness, the NPHC algorithm suffers from a large computational overhead as it employs numerous continuation traces to determine all solutions.

Past contributions are mostly focused on the study of only solvability space, i.e., without considering technological constraints. This manuscript aims to compute the power flow feasibility boundary with computational tractability, speed, and scalability to large networks. The key contributions of this paper are as follow:
\begin{enumerate}[I.]
    \item A general mathematical structure is proposed to extend the problem from solvability space towards the boundary of feasibility space by a set of complementary equations representing the operational inequalities.
    \item Compared to the eigenvector based transversality condition, a scalar transversality condition is used to characterize points on the solution space boundary. It results in a small computational overhead, does not suffer from initialization issues, and provides better scalability for large networks.
    \item This paper introduces the use of a Homotopy Continuation Method (HCM) referred to as ``Spherical Continuation'', which allows tracing the solution curves with sharp edges or non-convexities. An adaptive sphere size approach further improves the robustness and speed of this algorithm.
\end{enumerate}

\iffalse

%We comment this paragraph in the first submission due to space restrictions. 

The manuscript is organized in the following manner. Section.II describes the problem statement, followed by the mathematical structure of the power flow solution boundaries in Section.III. Details about proposed continuation with adaptive sphere strategy are provided in Section.IV. Whereas Section.V and VI provide detailed implementation and numerical results, respectively. Finally, conclusions and some future extensions of this work are presented in Section.VII.

\fi

% The very first letter is a 2 line initial drop letter followed
% by the rest of the first word in caps.
% 
% form to use if the first word consists of a single letter:
% \IEEEPARstart{A}{demo} file is ....
% 
% form to use if you need the single drop letter followed by
% normal text (unknown if ever used by the IEEE):
% \IEEEPARstart{A}{}demo file is ....
% 
% Some journals put the first two words in caps:
% \IEEEPARstart{T}{his demo} file is ....
% 
% Here we have the typical use of a "T" for an initial drop letter
% and "HIS" in caps to complete the first word.
%\IEEEPARstart{T}{his} demo file is intended to serve as a ``starter file'' for IEEE Communications Society journal papers produced under \LaTeX\ using
%IEEEtran.cls version 1.8b and later.
% You must have at least 2 lines in the paragraph with the drop letter
% (should never be an issue)
%I wish you the best of success.

%\hfill mds
 
%\hfill August 26, 2015

\section{Problem Statement}\label{math}
%Figure.1 provides a graphical depiction of the steady-state security regions in the space of $\lambda_{1}-\lambda_{2}$ (varying parameters, like power injections, voltage setpoints, etc.) to understand the significance of calculating boundaries of the power flow solution space in conjunction with power system security assessment.

\begin{figure}[t]
    \centering
    \includegraphics[scale=0.07]{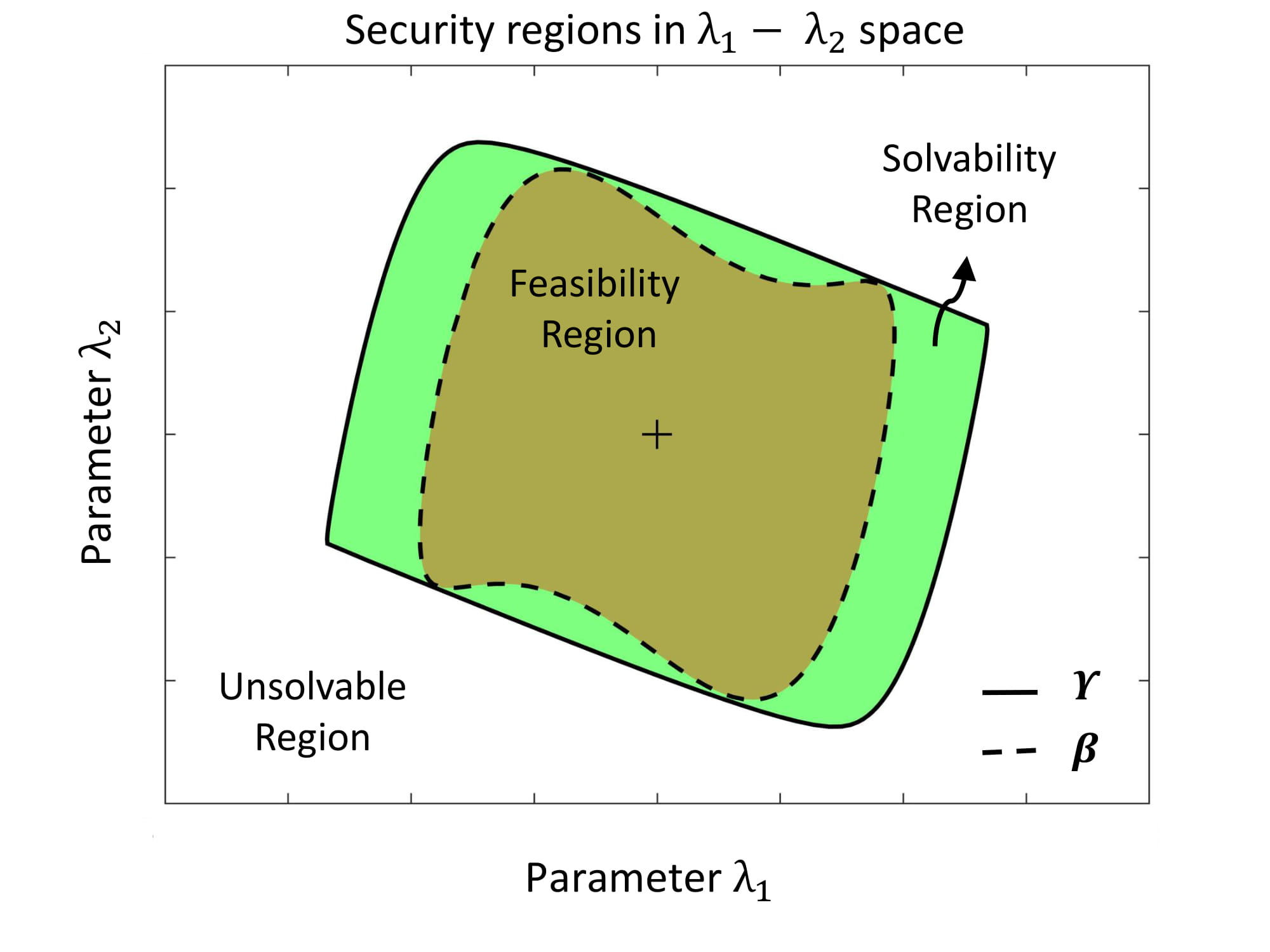}
    \caption{Steady-state security regions in $\lambda_{1}-\lambda_{2}$ space (varying parameters, like power injections, voltage setpoints, etc.).} 
    %in which the operator determines security, maneuverability or control actions needed to bring the operating point in the feasible domain.}
    \label{fig:feasbw}
\end{figure}

Here, Figure.\ref{fig:feasbw} depicts steady-state security regions such as feasibility, solvability, and unsolvable regions to illustrate the importance of the power flow solution map's boundedness in conjunction with the power system security assessment. %, we draw three different steady-state security regions in Figure.\ref{fig:feasbw}. 

First, the feasibility region describes the set of operating points for which the real-valued solution to power flow equations exists. The operational constraints like limits on voltage setpoints, generator power outputs, etc., are also satisfied. Usually, this is the optimal operating region for a grid. It is enclosed by a curve ``$\boldsymbol{\beta}$'' referred to as the voltage feasibility or just feasibility boundary. Secondly, the solvability region contains operating points for which the real-valued solution to the power flow equations still exists, but one or more operational constraints may be violated. The operation is still possible in this domain, but tractable measures must bring the operating point back to the feasible space. This region's boundary is outlined by ``$\boldsymbol{\gamma}$'' curve and generally known as the solvability boundary.  Finally, the area outside the solvability space is an unsolvable region, as the real-valued solution to the power flow equations disappears. Any operation in this space will result in system instability or voltage collapse.  

The quantification of the power flow solution boundaries (i.e., $\boldsymbol{\gamma}$ and $\boldsymbol{\beta}$ curves) can help Transmission System Operator (TSO) in security assessment by; i) determining the relative degree of security for a given operating state. ii) providing acceptable maneuverability for a continuously moving operating regime, iii) help decide optimum control actions to bring the operational state back to the feasible domain. %Information about solution boundaries can also be utilized in several other applications, i.e., contingency analysis, emergency controls, and maximizing the grid's transfer capabilities.

\subsection{Mathematical Formulation}
We can describe the power flow problem  with the following compact notation,  
\begin{equation}\label{eq:main}
    \mathcal{F} (\mathbf{x},\lambda ) = \mathbf{0}
\end{equation}
% \begin{equation}\label{eq:main}
%      f_{i}(x,\lambda)=0, \quad \quad i=1,...n.
% \end{equation}
Here, $\mathcal{F}(\mathbf{x},\lambda): \mathbb{R}^{n} \times \mathbb{R}^{p}\to \mathbb{R}^{n}$ constitute ``$n$'' number of nonlinear algebraic equations including both the power flows and the technological constraints, $\mathbf{x}\in\mathbb{R}^n$ is a vector of system variables (like voltage magnitudes and phase angles), while $\lambda \in\mathbb{R}^{p}$ are the system parameters (like limits on the nodal power injections, or limits on the line currents etc.).

\begin{figure}[!]
    \centering
    \includegraphics[scale=0.07]{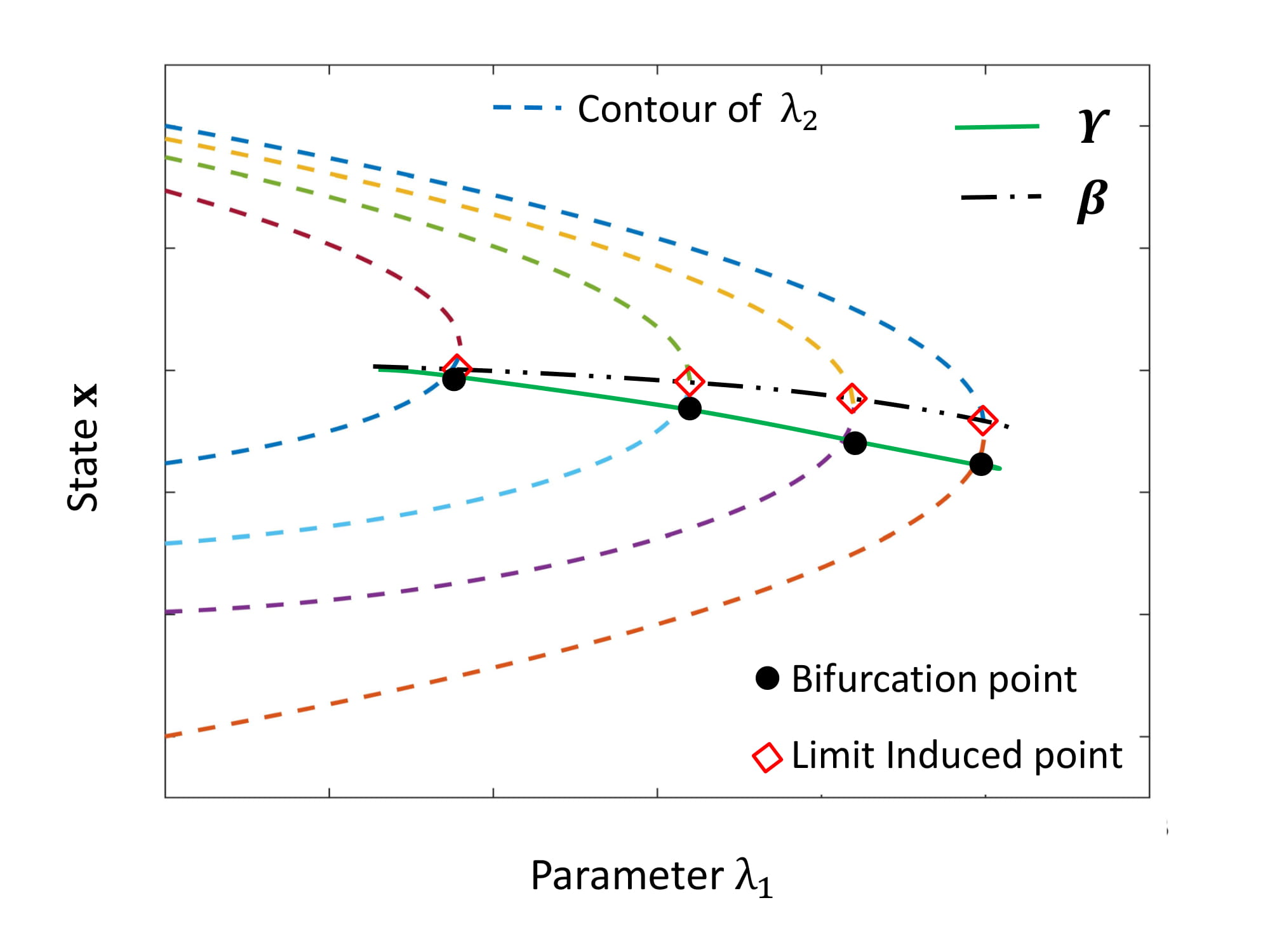}
    \caption{Power flow solution space boundaries in $\lambda_{1}-\lambda_{2}$ space.}
    \label{fig:ddddd}
\end{figure}

\paragraph{Solution space} The system in \eqref{eq:main} describes $n$-number of equations in $(n+p)$ variables. For ``$p$'' free varying parameters the solution space  of \eqref{eq:main} defines a $p-$dimensional manifold. Let $p=1$, solution of \eqref{eq:main} delineate $1-$dimensional curves. For example, dashed curves in Fig.\ref{fig:ddddd} represents solution of system in \eqref{eq:main} by releasing a single parameter $\lambda_{1}$, i.e., injected power while observing the state variable $\mathbf{x}$ like voltage magnitudes. These curves are analogous to PV or ``nose curves''  of the power flow problem. If another parameter $\lambda_{2}$ is also free to vary then the solution space of \eqref{eq:main} is a $2-$dimensional surface, dashed curves in Fig.\ref{fig:ddddd} are contours of such surface.  

\paragraph{Solution space boundary} This work aims to construct the boundary of the solution space of \eqref{eq:main}. Such boundary defines the boundedness of the power flow solution space. At any point on this boundary, the real-valued solution for the system in \eqref{eq:main} disappears either due to the unsolvability of power flow equations or due to the violation of some operational constraints. Mathematically, this condition can be enforced as, 
\begin{equation}\label{eq:sing}
    g(\mathbf{x}) = \det~ \nabla_\mathbf{x}\mathcal{F} (\mathbf{x},\lambda) = 0.
\end{equation}
Here \eqref{eq:sing} is referred to as transversality condition \cite{ali2017transversality} corresponding to the degeneracy of the $n \times n$ power flow Jacobain matrix $\nabla_\mathbf{x}\mathcal{F}$. There are several ways to enforce \eqref{eq:sing}, more details are given in the subsequent section. 

Introducing (\ref{eq:sing}) along with \eqref{eq:main}, makes the number of equations ($n+1$) in ($n+p$) variables. Following the same generalization, for ``$p$'' free varying parameters the solution space boundary of \eqref{eq:main} will be a ($p-1$)-dimensional manifold  in the parameter space of $\lambda$'s. If we consider $p=1$ then we have $(n+1)$ equations in $(n+1)$ variables and solution of \eqref{eq:main} together with \eqref{eq:sing} is just a single point on the solution space boundary.  In Fig.\ref{fig:ddddd} these boundary points are shown on both solvability boundary curve $\boldsymbol{\gamma}$ and also on the feasibility boundary $\boldsymbol{\beta}$. To sum up, the problem formulation for calculating the boundary of power flow solution space is described by solving the set of equations in \eqref{eq:main} along with condition \eqref{eq:sing}. The true solution boundary is a multi-dimensional hypersurface. Therefore, it is computationally infeasible to construct even for reasonably small networks. However, calculating a 1-dimensional slice of this high-dimensional hypersurface is computationally tractable and does provide a reliable way to assess the grid's security. The context of this work examines the structure of the solution boundary with $p=2$ resulting in a 1-dimensional boundary curves. The examples in Section.\ref{sec:numerics} also show a 2-dimensional solution boundary surface (i.e., for $p=3$).  
% \begin{subequations}\label{transequations}
% \begin{align}
%   \mathcal{F} (\mathbf{x},\lambda ) & = \mathbf{0},\\
%  g(\mathbf{x}) &= 0. 
% \end{align}
% \end{subequations}
%In the  subsequent sections we describe the system of equations and variables considered in $\mathcal{F} (\mathbf{x},\lambda ) = \mathbf{0}.$ 

%\section{System of Equations and Variables}\label{sec:eqvar}
%This section provides details about the system of equations and variables considered in $\mathcal{F} (\mathbf{x},\lambda ) = \mathbf{0}.$ 
\subsection{System of Equations and Variables}\label{powerflow}
In this section we describe the system of equations and variables considered in $\mathcal{F} (\mathbf{x},\lambda ) = \mathbf{0}.$ 
%The power flow model can be represented in both polar and rectangular coordinates. 
\subsubsection{Modeling Power Flows} Here, a rectangular formulation is adopted for the power flow model. %, as it simplifies the system of equations by neglecting the higher-order terms in the Taylor series \cite{tate2005comparison}. It also allows us to formulate inequality constraints into a polynomial form suitable to consider in \eqref{eq:main}. 
Let us consider the real representation of the power flow equations in terms of the real and imaginary part of the complex voltage phasor $\hat{V}_{i}\in \mathbb{C}^{n}$ as $\hat{V}_{i}=V^{r}_{i}+\boldsymbol{j} V^{m}_{i}$. Let $n$ be the total number of buses in the network with $\mathcal{N}=\{1,2,\dots,n\}$ representing the set of all buses, while $\mathcal{L}$ is the set of load buses (PQ), and $\mathcal{G}$ is the set of generator buses (PV). For each bus, $i\in \mathcal{N}$ except the slack bus, the active power balance can be stated as, 
\begin{equation}
  \begin{aligned}\label{eq:pi}
%P_{i} (x)=
\sum_{k}%^{n} 
\Big\{V^{r}_{i}(G_{ik}V^{r}_{k}-B_{ik}V^{m}_{k})
+ V^{m}_{i}(G_{ik}V^{m}_{k}+B_{ik}V^{r}_{k}) \Big\}-\\
-\lambda(P^{\mathrm {gen}}_{i}-P^{\mathrm{load}}_{i}) =0.
  \end{aligned}
\end{equation}
Similarly, the reactive power balance at each load bus $i\in \mathcal{L}$ can be defined as,
\begin{equation}
  \begin{aligned}\label{eq:qi}
%Q_{i} (x)=
\sum_{k}%^{n} 
\Big\{V^{m}_{i}(G_{ik}V^{r}_{k}-B_{ik}V^{m}_{k}) -V^{r}_{i}(G_{ik}V^{m}_{k}+B_{ik}V^{r}_{k}) \Big\}-\\
-\lambda(Q^{\mathrm {gen}}_{i}-Q^{\mathrm{load}}_{i}) =0.
  \end{aligned}
\end{equation}
The parameter $\lambda$ represents the boundness of the nodal power injections, while subscripts ``gen'' and ``load'' express bus loading and generation levels. In  \eqref{eq:pi} and \eqref{eq:qi} $G_{ij}$ and  $B_{ij}$ denote the real and imaginary entries of the complex admittance matrix $\hat{Y}_{ij}=G_{ij}+\boldsymbol{j} B_{ij}$. For each generator, $i\in \mathcal{G}$ the fixed voltage levels are characterized by an additional equation,
\begin{equation} \label{eq:vi}
(V^{r}_{i})^{2}+(V^{m}_{i})^{2}-|\hat{V_{i}}|_{\mathrm{ref}}^{2}=0. 
\end{equation}
In \eqref{eq:vi} $|\hat{V_{i}}|_{\mathrm{ref}}$ is the specified voltage magnitude at the $i^{th}$ generator bus.

\subsubsection{Voltage magnitude limits}\label{feasbi}
If the system in \eqref{eq:main} considers only the power flow equations i.e., \eqref{eq:pi}, \eqref{eq:qi} and \eqref{eq:vi}, then the problem formulation is limited to the boundary of the solvability space. Typically, a real power system is operated with technological constraints. Usually, these constraints are given in the form of inequalities with maximum and minimum bounds. To extend the problem to the boundary of the feasibility space (i.e., $\boldsymbol{\beta}$ curve from Fig.1\ref{fig:feasbw}), we reformulate inequality constraints into a set of polynomial equations suitable to consider in $\eqref{eq:main}$.
%\paragraph{Voltage magnitude limits}\label{vglimit} 
For example, at each load bus $i\in \mathcal{L}$, voltage magnitudes are given within specified bounds as,
%This constraint restricts the solution of \eqref{eq:pi} such that voltage magnitude at each PQ bus remains within certain bound as, 
\begin{equation}\label{eq:inqvpq}
    |\hat{V}_{i}|_{\mathrm{min}} \leq\ |V_{i}|_{\mathrm{cal}}\leq |\hat{V}_{i}|_{\mathrm{max}}
\end{equation}
Here $|V_{i}|_{\mathrm{cal}}$ specifies the calculated voltage magnitude level. While $|\hat{V}_{i}|_{\mathrm {max}}$ and $|\hat{V}_{i}|_{\mathrm {min}}$ represent maximum and minimum bound on the voltage magnitude. To have a solvable solution of \eqref{eq:inqvpq}, slack variables were introduced to convert this inequality constraint into two different equality equations corresponding to  $|\hat{V}_{i}|_{\mathrm{max}}$ and $|\hat{V}_{i}|_{\mathrm{min}}$ respectively,
\begin{subequations}
\label{eq:vlimits}
\begin{align}
(V_{i}^{r})^{2}+(V_{i}^{m})^{2}-|\hat{V}_{i}|_{\mathrm{max}}^{2}+(\overline{s}_{i})^{2} &=0   \label{eq:sla11}\\
(V_{i}^{r})^{2}+(V_{i}^{m})^{2}-|\hat{V}_{i}|_{\mathrm{min}}^{2}-(\underline{s}_{i})^{2} &=0 \label{eq:sla12}
\end{align}
\end{subequations}
In \eqref{eq:vlimits}, $\overline{s}_{i}$ and $\underline{s}_{i}$ denote the slack variables related to each limit. From \eqref{eq:vlimits}, one can conclude that the real solution for slack variables $\overline{s}_{i}$ and $\underline{s}_{i}$ exits only if voltage levels lie within acceptable limits, and vice-versa, such a solution does not exist if the limits are violated. Thus, we have reduced the feasibility problem to the solvability problem of a certain equivalent system.  

Here, we consider the constraints of the type \eqref{eq:inqvpq}, however any other constraints like limits on generator power outputs or transmission line thermal limits can be reformulated using the same slack variable methodology \cite{ali2019fast} similarly to \eqref{eq:vlimits}. It is worth noticing that even with additional equations, the given methodology preserves the power flow model's sparsity characteristics. To summarize, the system in \eqref{eq:main} considers equations, \eqref{eq:pi}, \eqref{eq:qi}, \eqref{eq:sla11}, and \eqref{eq:sla12} for every PQ bus, while two equations - \eqref{eq:qi} and \eqref{eq:vi} - are defined for every PV bus. The variable vector $\mathbf{x}$ contains $n$ components, namely variables $V_i^r , V_i^m , \overline{s}_i , \underline{s}_i$ for every PQ bus and variables $V_i^r , V_i^m$ for every PV bus.

\section{Proposed Methodology}
In the present section, we first outline the choice of transversality condition $g(\mathbf{x})$ followed by the discussion of the proposed algorithm. 
\subsection{Transversality condition $g(x)$}\label{g(x)tran}
Normally, the transversality condition is applied using the determinant as in \eqref{eq:sing}. It is simple to implement but suffers from a high computational burden, numerical stability issues, and poor scalability to larger networks. The conventional way to enforce this condition uses a right eigenvector corresponding to the null space of $ \nabla_\mathbf{x}\mathcal{F} (\mathbf{x},\lambda)$ as,
\begin{equation}\label{eq:eigenvector}
g_{\mathrm{eig}}(\mathbf{x}, \mathbf{y}) =
\begin{bmatrix}
  \nabla_\mathbf{x}\mathcal{F} (\mathbf{x},\lambda) \cdot \mathbf{y} \\
  \mathbf{y}^\top \mathbf{y} - 1
\end{bmatrix} = \mathbf{0}. 
\end{equation}
In \eqref{eq:eigenvector} $\mathbf{y}$ is a normalised eigenvector corresponding to zero eigenvalue of $\nabla_\mathbf{x}\mathcal{F} (\mathbf{x},\lambda)$. Introducing $g(\mathbf{x})$ through \eqref{eq:eigenvector} has the benefit of preserving the power flow problem's sparsity structure. But, the convergence is sensitive to a good initial guess of the eigenvector $\mathbf{y}$ and displays poor scalability as dimension of system variables increases by a factor of two.

Several choices of $g(\mathbf{x})$ are explored in \cite{ali2017transversality} with fast algorithms to evaluate their gradients. We employ a scalar transversality condition based on the Singular Value Decomposition (SVD) due to the small computational cost, numerical stability, and scalability to large networks. Let the SVD of the $n\times n$ Jacobian matrix $\boldsymbol{J}$ i.e., $\nabla_\mathbf{x}\mathcal{F} (\mathbf{x},\lambda)$ be,
\begin{equation}\label{eq4}
\boldsymbol{J} = \mathbf{U} \boldsymbol{\Sigma} \mathbf{{V}}^{\top}= \sum_{i=1}^{n} {\sigma_{i}}\mathbf{u}_{i}\mathbf{v}_{i}^{\top}.
\end{equation}
Here $\mathbf{U}$ and $\mathbf{V}$ are the orthogonal matrices with same dimensions as $\boldsymbol{J}$, while $\boldsymbol{\Sigma}$ is a diagonal matrix consisting of singular values such that $\sigma_{i} \geq 0,~i=1,2,\cdots,n$. At the solution boundary the Jacobin becomes singular, thus the minimum singular value of $\boldsymbol{J}$ becomes equal to zero i.e.,~$\sigma_{nn}=0$. Consequently, one can enforce the following necessary condition for characterizing the solution of \eqref{eq:main} on the feasibility boundary,
\begin{equation}\label{eq:svd}
    g_{\mathrm{svd}}(\mathbf{x})= \mathbf{u}^{\top}_{n}\boldsymbol{J}\mathbf{v}_{n} = 0.
\end{equation}
Here, $\mathbf{u}_{n}$ and $\mathbf{v}_{n}$ correspond to the $n^{th}$ left and right singular vectors. Introducing $g(\mathbf{x})$ by \eqref{eq:svd} has following advantages over \eqref{eq:eigenvector}; i) small computational overhead as only one additional variable is introduced, ii) does not require any guess for initialization compared to $g_{\mathrm{eig}}$, iii) the singular values are less sensitive to numerical perturbations. iv) finally, better scalability as this condition requires computing SVD only for the lowest singular value, which can be performed with a breakneck speed for sparse structured matrices like Jacobian $\boldsymbol{J}$. One can now describe the solution of system in \eqref{eq:main} on the boundary by enforcing condition \eqref{eq:svd}. 

% \[\boldsymbol{\beta} = \Big\{(\mathbf{x},\lambda): \forall~\mathbf{x}\in \mathbb{R}^{n},~\lambda \in \mathbb{R}^{p}, \mathcal{F}(\mathbf{x},\lambda) = 0,~ g_\mathrm{svd} (\mathbf{x}) = 0 \Big\}\]

\subsection{Homotopy Continuation Method}\label{overview}
The general set of equations defining the points on the boundary of the power flow solution space can be stated as,
\begin{equation}\label{eq:h(z)1}
    \mathcal{H}(\mathbf{z}) = \mathbf{0}.
\end{equation}
Where,
 \begin{equation}\label{eq:dddsaq}
         \mathcal{H} (\mathbf{z})=\begin{bmatrix}
         \mathcal{F} (\mathbf{x},\lambda )\\
         \\
         g_\mathrm{svd} (\mathbf{x})
         \end{bmatrix}\in\mathbb{R}^{n+1}, \quad
    \mathbf{z} = \begin{bmatrix}
            \mathbf{x}\\
            \\
           \lambda\\
         \end{bmatrix}\in\mathbb{R}^{n+p}
\end{equation}
Using the same generalization, if there are ``$p$'' free varying parameters i.e., $\lambda \in \mathbb{R}^{p}$, the solution of \eqref{eq:h(z)1} corresponds to a $(p-1)-$dimensional manifold. This paper considers formulation with $p=2$ resulting in a curve or 1-dimensional curve of solution points. 
%in the $(n+2)-$dimensional Euclidean space. %The set of all $z = (x_{1},x_{2},\dots,\lambda)^{\top}$ that satisfy \eqref{eq:h(z)1} will generally consist of one or more solution curves in the $(n+2)-$dimensional Euclidean space. 

Generally, Homotopy Continuation Methods (HCMs) are deployed for tracing solution curves, with some suitable path-tracking algorithms. Namely, i) Euler Homotopy or Predictor-Corrector (PC) algorithms, and ii) Simplicial Continuation or Piecewise-Linear (PL) algorithms \cite{yamamura1993simple, allgower2003introduction}. %The main difference is the condition on the homotopy map's smoothness defined by the set of equations in \eqref{eq:h(z)1}. The PC algorithms require a smooth map of \eqref{eq:h(z)1}, and vice versa, such a condition is not necessary for PL algorithms. 
Here, we describe a PC algorithm to trace the homotopy path, as the underlying equations in \eqref{eq:h(z)1} represent a smooth map. The PC algorithms are also more efficient in higher-dimensional cases compared to the PL algorithms. Conventionally, the PC algorithms apply tangents as predictors and hyperplanes as correctors (like in \cite{Hiskens:2001ie}). Such an algorithmic routine suffers from difficult implementation, convergence issues, and manifest slow performance in the vicinity of sharp turning points or non-convex sections of the solution curve \cite{yamamura1993simple}. This paper implements an adaptive Spherical Continuation (SC) algorithm for curve tracing using spheres, with a precise geometrical description, ease of implementation, and computational tractability. 
%\subsection{Predictor-Corrector Procedure}

\subsubsection{Spherical Continuation Path Tracking}
\begin{figure}
    \centering
    \includegraphics[scale=0.07]{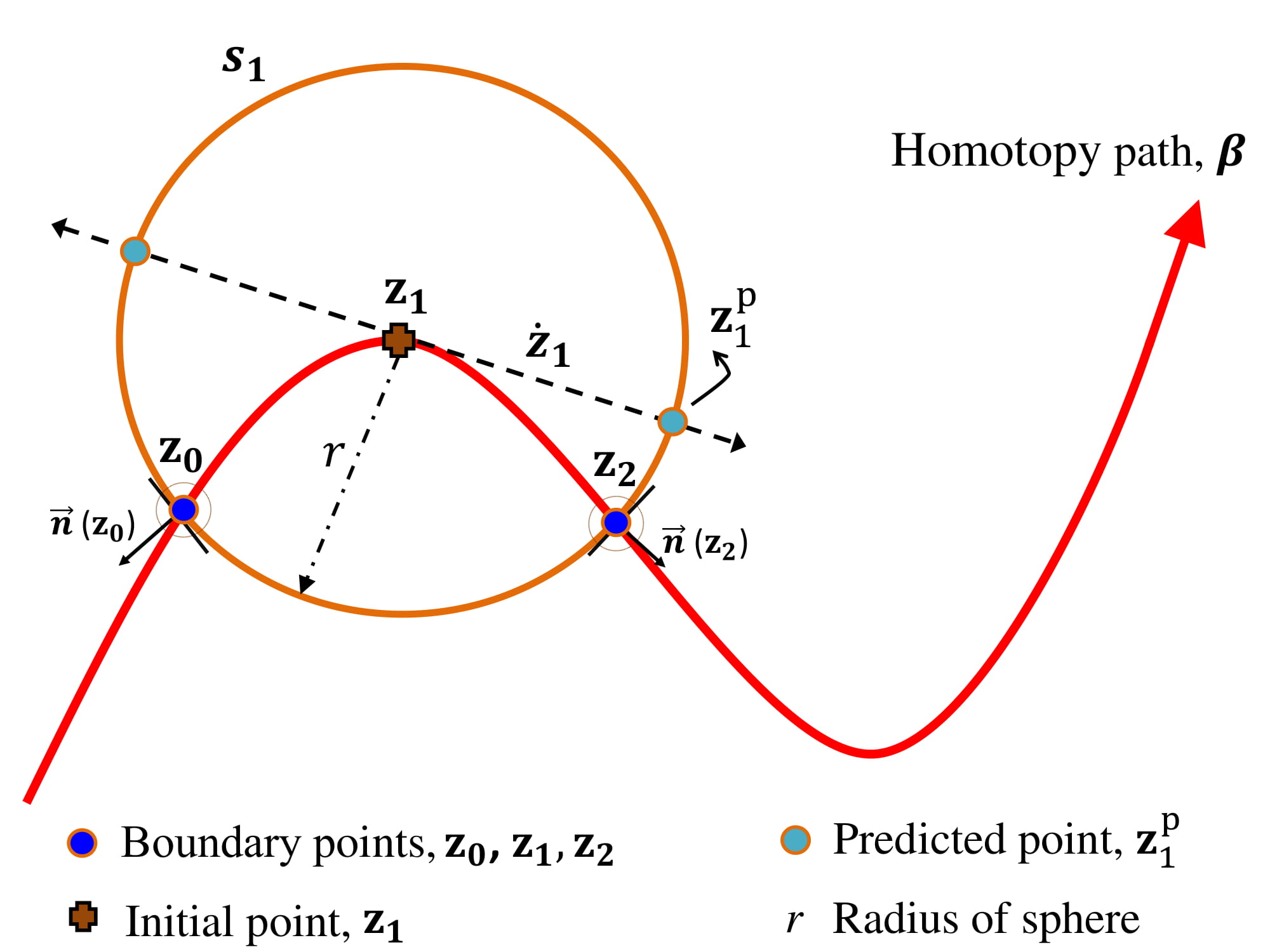}
    \caption{Intersection of sphere $\boldsymbol{s_{1}}$ with a homotopy path $\boldsymbol{\beta}$.}
    \label{Sphere_1}
\end{figure}

\begin{figure}
    \centering
    \includegraphics[scale=0.07]{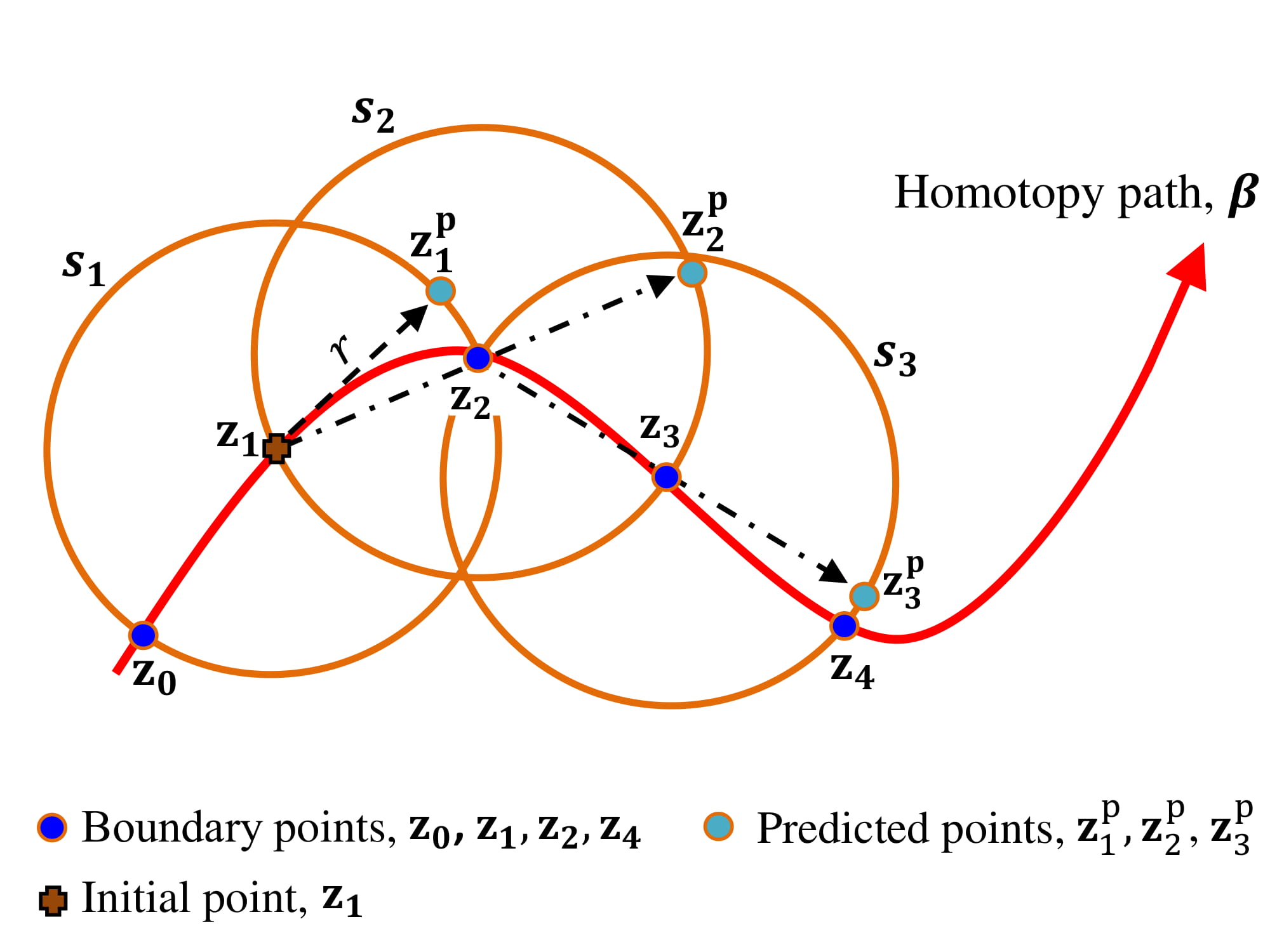}
    \caption{Spherical continuation, path tracking procedure.}
    \label{fig:sphere}
\end{figure}

Let us consider a sphere ``$\boldsymbol{s_{1}}$'' of radius ``$r$'' as shown in Fig.\ref{Sphere_1}, centered at $\mathbf{z}_{1}$ on the homotopy path ``$\boldsymbol{\beta}$'' defined by \eqref{eq:h(z)1}.  Its outline must intersect the homotopy curve at least at two possible points, one in the forward direction ``$\mathbf{z}_{0}$'', and another one in the backward direction ``$\mathbf{z}_{2}$''. The given formulation of the homotopy system \eqref{eq:h(z)1} has $(n+1)$ equations in $(n + 2)$ variables, it means that we should include an additional equation representing a $(n+2)$-dimensional hypersphere to trace the homotopy path. Thus, the equation of the hypersphere intersecting the homotopy path $\boldsymbol{\beta}$ can be introduced as,
\begin{equation}\label{hyperspherell}
    \mathcal{S} (\mathbf{z},r) = \sum_{m}^{n+2} (\mathbf{z}^{m} - \mathbf{z}^{m}_{1})^{2} - r^{2} = 0.
\end{equation}
Figure \eqref{fig:sphere} shows the Spherical path tracking of the homotopy curve $\boldsymbol{\beta}$. To trace the succeeding points, i.e., $\mathbf{z}_2, \mathbf{z}_3, \cdots$ on the homotopy curve, the PC algorithm starts from an Euler Predictor followed by a Generalized Newton-Raphson Corrector. 

\paragraph{Euler predictor} The algorithm is driven from an initial point $\mathbf{z}_{1}$ on the $\boldsymbol{\beta}$ curve that satisfies $\mathcal{H}(\mathbf{z}_{1})=0$. To compute the next point $\mathbf{z}_2$ on the curve, first a predictor point $\mathbf{z}_{1}^{p}$ is calculated on the sphere $\boldsymbol{s_{1}}$ through an Euler predictor. 

Let us assume an arc-length parametrization ``$\alpha$'' of the homotopy curve $\boldsymbol{\beta}$ is a smooth function $\mathbf{z} (\alpha): \mathbb{R}^{+} \to \boldsymbol{\beta}$ such that $\mathbf{z}(0) = \mathbf{z}_1$,  and $\mathcal{H} \big(\mathbf{z(\alpha)}\big) = 0$.  The tangent vector $\dot{\mathbf{z}} (\alpha)$ always has a unit length for all $\alpha \in \mathbb{R}^{+}$. The aforementioned conditions for $\mathbf{z}(\alpha)$ make it a solution of the following system of differential equations,
\begin{subequations}\label{eq:eeet}
\begin{align}
    \nabla_\mathbf{z}\mathcal{H (\mathbf{z} (\alpha))}  \cdot \dot{\mathbf{z}} (\alpha) & = \mathbf{0}\\
    \dot{\mathbf{z}} (\alpha)^{\top} \cdot \dot{\mathbf{z}} (\alpha) & =1.\\
    \mathbf{z}(0) & = \mathbf{z}_1
\end{align}
\end{subequations}
Where $ \nabla_\mathbf{z}\mathcal{H (\mathbf{z} (\alpha))}$ is the ${(n+1)\times (n+2)}$ Jacobian matrix and has a rank $(n+1)$ by assumption. While $\dot{\mathbf{z}} (\alpha)$ is the normalized tangent vector spanning the one-dimensional kernel $\mathrm{ker} \Big(  \nabla_\mathbf{z}\mathcal{H (\mathbf{z} (\alpha))}   \Big)$. From \eqref{eq:eeet} it is clear that the kernel of the Jacobian $ \nabla_\mathbf{z}\mathcal{H (\mathbf{z} (\alpha))}$ has two vectors of unit norm which correspond to the two possible directions of traversing the $\boldsymbol{\beta}$ curve (see Fig.\ref{Sphere_1}). To specify the direction in which the curve is traversed following condition is introduced along with \eqref{eq:eeet}
\begin{equation}\label{eq:signd}
    \mathrm{det}~\begin{pmatrix}
       \nabla_\mathbf{z}\mathcal{H (\mathbf{z} (\alpha))}\\
       \dot{\mathbf{z}} (\alpha)^{\top}
    \end{pmatrix} > 0
\end{equation}
% {\color{blue}The tangent vector $\dot{\mathbf{z}} (\alpha)$ can be calculated numerically via QR decomposition if,
% \begin{equation}
%      \Big( \nabla_\mathbf{z}\mathcal{H} (\mathbf{z}) \Big)^{\top} \big|_{\mathbf{z} = \mathbf{z}_{1}} = Q\begin{pmatrix}
%      R\\
%      \mathbf{0}
%     \end{pmatrix} 
% \end{equation}
% Here $Q$ is an orthogonal matrix composed of $(n+2)$ orthonormal vectors $Q=[\mathbf{q}_1,\mathbf{q}_2,\cdots ,\mathbf{q}_{n+2}]$, and $R$ is an upper-triangular matrix with full rank. Let $\delta \mathbf{z}$ is the $n^{th}$ orthonormal vector such that $\delta \mathbf{z} = \mathbf{q}_{n+2} = Q\cdot \boldsymbol{e}_{n+2}$  (with $\boldsymbol{e}_{n+2} = (0,0,\dots,0,1)^{\top}\in \mathbb{R}^{(n+2)}$), then it satisfy the following,}
The tangent vector $\dot{\mathbf{z}} (\alpha)$ can be calculated numerically via SVD decomposition if,
\begin{equation}\label{eq:tansvd}
     \left.\Big( \nabla_\mathbf{z}\mathcal{H} (\mathbf{z}) \Big)^{\top} \right|_{\mathbf{z} = \mathbf{z}_{1}} = \mathbf{M}\begin{pmatrix}
     \boldsymbol{D}\\
     \mathbf{0}
    \end{pmatrix} \mathbf{N}^{\top}.
\end{equation}
Here $\mathbf{M}$ is a ${(n+2)\times (n+2)}$ left orthogonal matrix composed of $(n+2)$ orthonormal vectors $ \mathbf{M} =[\mathbf{m}_1,\mathbf{m}_2,\cdots ,\mathbf{m}_{n+2}]$, and $\boldsymbol{D}$ is a ${(n+1)\times (n+1)}$ full rank diagonal matrix consisting of $(n+1)$ singular values. While  $\mathbf{N}$ is a ${(n+1)\times (n+1)}$ right orthogonal matrix composed of $(n+1)$ orthonormal vectors $ \mathbf{N} =[\mathbf{n}_1,\mathbf{n}_2,\cdots ,\mathbf{n}_{n+1}]$. After some algebraic manipulation, the expression in \eqref{eq:tansvd} can be stated as,
\begin{equation}
     \nabla_\mathbf{z}\mathcal{H} (\mathbf{z}) \big|_{\mathbf{z} = \mathbf{z}_{1}} \mathbf{M}  = \mathbf{N} \Big(
     \mathbf{D} \quad \mathbf{0} \Big).
\end{equation}
Let $\delta \mathbf{z}$ is the $n^{th}$ orthonormal vector such that $\delta \mathbf{z} = \mathbf{m}_{n+2} = \mathbf{M}\cdot \mathbf{e}_{n+2}$  (with $\mathbf{e}_{n+2} = (0,\dots,0,1)^{\top}\in \mathbb{R}^{n+2}$) such that:
\begin{subequations}
\begin{align}
    \nabla_\mathbf{z}\mathcal{H} (\mathbf{z}) \big|_{\mathbf{z} = \mathbf{z}_{1}}\cdot \delta \mathbf{z} &= \mathbf{0} \\
    \delta \mathbf{z}^{\top} \cdot  \delta \mathbf{z} & = 1.
\end{align}
\end{subequations}
% {\color{blue} The condition \eqref{eq:signd} for traversing the curve in a forward direction can be interpreted by the sign of the $\Big( \mathrm{det}~(Q)\cdot \mathrm{det}~(R) \Big)$, which can be computed without any additional computational burden.}
Condition \eqref{eq:signd} for traversing the curve in a forward direction can be interpreted by the sign of the $\Big( \mathrm{det}~(\mathbf{M})\cdot \mathrm{det}~(\mathbf{D}) \cdot \mathrm{det}~(\mathbf{N}) \Big)$, which can be computed without any additional computational cost. Finally, the predicted point on the sphere $\boldsymbol{s}_{1}$ with radius $r$ is given by,
\begin{equation}\label{eq:predictorx}
    \mathbf{z}_{1}^{p} = \mathbf{z}_{1} + r\cdot \delta \mathbf{z}
\end{equation}
\paragraph{Generalized Newton Corrector} The predicted solution $\mathbf{z}_{1}^{p}$ is just an estimate of the actual point on the curve. Hence, a correction step is essential to calculate an exact point on the curve. The subsequent set of equations are solved in the corrector step to calculate a forward solution $\mathbf{z}_{2}$ on the homotopy path, 
\begin{subequations}\label{eq:hypersphere}
\begin{align}
    \mathcal{H} (\mathbf{z}) &= \mathbf{0},\\
\mathcal{S} (\mathbf{z}, r) = \sum_{m}^{n+2} (\mathbf{z}^{m} - \mathbf{z}^{m}_{1})^{2} - r^{2} &= 0.    
\end{align}
\end{subequations}
% \begin{equation}
% \left\lVert {\mathbf{z} - \mathbf{z}_{1}} \right\rVert^{2}
% \end{equation}
The system in \eqref{eq:hypersphere} forms a  set of $(n+2)$ equations in $(n+2)$ unknowns, which can be solved using any standard technique, such as Newton-Raphson. Linearizing \eqref{eq:hypersphere} at $\mathbf{z}_1^{p}$ by first order Taylor series expansion leads to,
\begin{equation}\label{eq:newtona}
    \Delta \mathbf{z} = \mathbf{z} - \mathbf{z}_1^{p} = \left[
    \begin{pmatrix}
      \nabla_\mathbf{z} \mathcal{H} (\mathbf{z})\\
      \nabla_\mathbf{z} \mathcal{S} (\mathbf{z},r)^{\top}\\
    \end{pmatrix}^{-1}
    \begin{bmatrix}
      \mathcal{H} (\mathbf{z})\\
      \mathcal{S} (\mathbf{z},r)
    \end{bmatrix}\right]_{\mathbf{z}=\mathbf{z}_1^p}
\end{equation}
The Newton point $\mathcal{N} (\mathbf{z})$ for approximating the solution of \eqref{eq:newtona} is given by,
\begin{equation}
    \mathcal{N} (\mathbf{z}) :=\mathbf{z}_1^{p} + \eta \left[
    \begin{pmatrix}
      \nabla_\mathbf{z} \mathcal{H} (\mathbf{z})\\
      \nabla_\mathbf{z} \mathcal{S} (\mathbf{z},r)^{\top}\\
    \end{pmatrix}^{-1}
    \begin{bmatrix}
      \mathcal{H} (\mathbf{z})\\
      \mathcal{S} (\mathbf{z},r)
    \end{bmatrix}\right]_{\mathbf{z}=\mathbf{z}_1^p}
\end{equation}
The ``$\eta$'' is Newton step size, which should be chosen to be small enough to ensure global convergence. The map $\mathcal{N} (\mathbf{z}) : \mathbb{R}^{(n+2)} \to \mathbb{R}^{(n+2)}$ is referred to as Newton map, that is iteratively applied to $\mathbf{z}$ until a certain convergence criterion 
\begin{equation}\label{eq:Newtonmap}
    \|  \mathcal{N}^{\kappa}(\mathbf{z}) -    \mathcal{N}^{\kappa-1}(\mathbf{z}) \| \leq \epsilon, \quad \kappa = 1,2,\cdots 
\end{equation}
is met. In \eqref{eq:Newtonmap} $\epsilon$ denotes the desired tolerance level. Usually, the Newton method converges quadratically with a predictor point $\mathbf{z}_{1}^{p}$ sufficiently close to $\mathbf{z}_2$. If iterations fail to converge, then a step size adjustment is required. 

%The map $\mathcal{N} : \mathbb{R}^{(n+2)} \to \mathbb{R}^{(n+2)}$ will also be called the Newton map. From the potentially inaccurate prediction $\mathbf{z}_{1}^{p}$ produced by Euler’s predictor, the map $\mathcal{N}$ is iteratively applied to $\mathbf{z}_{1}^{p}$ until a certain convergence criterion is met. As stated formally, the result of the correction is,
% \begin{equation}
%     \mathbf{z}  = \mathcal{N}^{k}(\mathbf{z}_{1}^{p}) = \mathcal{N} \circ \mathcal{N} \dots \circ  \mathcal{N}(\mathbf{z}_{1}^{p})
% \end{equation}
% Where $k \in \mathbb{Z}^{+}$ is determined by the convergence criterion, it is well known that this generalized Newton method converges quadratically for $\mathbf{z}_{1}^{p}$ sufficiently close to $\mathbf{z}$ under the smoothness assumption (see, e.g., [88] and [89]). Therefore, the shrinking distance $\|  \mathcal{N}^{j}(\mathbf{z}_{1}^{p}) -    \mathcal{N}^{j-1}(\mathbf{z}_{1}^{p}) \|$ between successive points produced by the iterations can be used as a criterion for convergence. Of course, if the iterations fail to converge, one must go back to adjust the step size for the Euler’s predictor.

%At each continuation step, the preceding solution is used as a center of the new sphere, as depicted in Fig.\ref{fig:sphere}, thus eliminating the requirement for an explicit tangent based predictor. It reduces the overall computational burden as the predicted point $\mathbf{z}_{i}^{p}$ at each continuation step ``$i$'' can be estimated by extrapolating between the consecutive spheres' centers. 

\subsubsection{Reversion Phenomenon} From Fig.\ref{Sphere_1}, it is clear that the sphere $\boldsymbol{s_{1}}$ intersects the homotopy curve at two possible points, one in the forward direction $\mathbf{z}_{0}$ and another in the backward direction $\mathbf{z}_{2}$. Therefore, there is still a chance that even with a predictor point ${\mathbf{z}}_1^p$ in the close vicinity of a forward solution $\mathbf{z}_{2}$, the Newton-Raphson in the corrector step can converge to a preceding solution for the system in \eqref{eq:hypersphere}. This phenomenon is known as the reversion effect or backtracking \cite{yamamura1993simple}.

The backtracking phenomenon can be detected by evaluating the normal vector's direction at each point of intersection. For the $\boldsymbol{s_{1}}$ sphere, the following expression allows us to evaluate the directions of the normal vector,
\begin{equation}\label{eq:reversal}
    \bigg( \cos^{-1}\Big(\frac{\nabla_\mathbf{z} \mathcal{S}(\mathbf{z},r)} {||\nabla_\mathbf{z} \mathcal{S}(\mathbf{z},r)||}\Big)  \bigg)_\mathbf{z_{0}} =   \bigg( \cos^{-1}\Big(\frac{\nabla_\mathbf{z} \mathcal{S}(\mathbf{z},r)} {||\nabla_\mathbf{z} \mathcal{S}(\mathbf{z},r)||} \Big) \bigg)_\mathbf{z_{2}}
\end{equation}
In \eqref{eq:reversal} $\nabla_\mathbf{z} \mathcal{S}(\mathbf{z},r) \in \mathbb{R}^{n+2}$ is the gradient of the corrector equation \eqref{hyperspherell}, while  $\cos^{-1}$ function is used to find the angle with respect to each variable in the normal vector. Once the reversion is detected, a slightly perturbed predicted solution $\tilde{\mathbf{z}}_1^p$ is used as an initial guess in the corrector step by increasing the sphere's radius by a small amount ``$\delta r$'', resulting in a  forward solution. 

\subsubsection{Adaptive sphere strategy}\label{adap}
Traditionally, HCMs are reprimanded for slowness in path tracking. Hence, an adaptive sphere size strategy is devised to improve the computational speed and modify the algorithm to trace the boundary curves with sharp edges. The proposed approach is based on understanding the radius of curvature of the homotopy path.  Larger the curvature radius flatter the curve. Contrarily, the smaller the radius of curvature steeper the curve. As the system in \eqref{eq:h(z)1} defines a homotopy curve by the intersection of $(n+1)$ implicit hypersurfaces, the radius of curvature formula can be stated. as\cite{goldman2005curvature}, 
\begin{equation}\label{eq:tau}
    \tau = \frac{\| \dot{\mathbf{z}} (\alpha) \|^{3}}{\Big| \Big( \dot{\mathbf{z}} (\alpha)^{\top} \cdot \nabla_\mathbf{z}~ \dot{\mathbf{z}} (\alpha) \Big)^{\top}\wedge \dot{\mathbf{z}} (\alpha) \Big|}
\end{equation}
Here $\dot{\mathbf{z}} (\alpha)$ is the normalized tangent vector calculated numerically through \eqref{eq:tansvd}, while $\nabla_\mathbf{z}~ \dot{\mathbf{z}} (\alpha)  \in \mathbb{R}^{(n+2)\times (n+2)}$ is the gradient of tangent vector. And operator ``$\wedge$'' is the wedge product. Here, we propose an asymptotic function of $\tau$ to calculate the radius of sphere $r$ from a maximum and minimum bound,
\begin{equation}\label{eq:radius}
    r = \mathrm{tanh} (\tau)= \frac {  \mathcal{A} \mathrm{e}^{\tau} - \mathcal{B} \mathrm{e}^{-\tau}}{\mathcal{C} \mathrm{e}^{\tau} + \mathcal{D} \mathrm{e}^{-\tau}}
\end{equation}

In \eqref{eq:radius} the values of $\mathcal{C}$ and $\mathcal{D}$ are 1, while $\mathcal{A}$ and $\mathcal{B}$ correlates to the maximum and minimum bound on the value of the sphere's radius, respectively.
% \begin{equation}
%     r = r_{0} \Big( 1- (1-\mathbf{\Gamma})~e^{-|\frac{\mu_{i}}{\mu_{i+1}}|} \Big)
% \end{equation}
It was noted that with the proposed methodology, the spherical path tracking scheme becomes faster and more adaptable for tracing the boundary curves with sharp turns.

\subsection{Algorithm Implementation}\label{implementation}
A two-stage implementation is adopted to trace the solution space boundary of \eqref{eq:main}. First, a single-point solution is calculated by varying only a single parameter, i.e., $\lambda_1$.  Hence \eqref{eq:h(z)1} becomes a set of $(n+1)$ equations in $(n+1)$ unknowns. To solve such a set of equations, one can use the Transversality Enforced Newton-Raphson (TENR) method from \cite{ali2017transversality}.  As the system in \eqref{eq:h(z)1} consists of non-linear algebraic equations, this means one or more isolated solutions exists. Each of these isolated solutions associates with a distinct solution boundary curve. Relaxed TERN iterations are deployed to identify these various roots, later to be used as an initial starting point for tracing the corresponding solution branch. 

In the next stage, the second parameter $\lambda_{2}$ is also free to vary; this results in a underdetermined system with $(n+1)$ equations in $(n+2)$ variables. Therefore, the solution of \eqref{eq:h(z)1} becomes a 1-dimensional curve, which is traced by the proposed SC method that uses the initial point, i.e., $\mathbf{z}_{i} $ from the first stage. Detailed procedure for the SC path tracking is presented in Algorithm \ref{euclid}.

\begin{algorithm}
  \caption{Spherical Continuation Path Tracking}\label{euclid}
  \begin{algorithmic}[1]
  \Procedure{Predictor-Corrector}{$\mathbf{z}_{i}, r$} %\Comment{The g.c.d. of a and b}
  \Statex \textbf{input}
  \State \textbf{begin}
  \State $z_{i}\in \mathbb{R}^{n+2},~\text{such that}~\mathcal{H}(\mathbf{z_{i}}) = 0;$ \Comment{initial point}
  \State $r>0;$ \Comment{initial  sphere's radius}
  \State \textbf{end};
  \Statex \textbf{repeat} \Comment{predictor step}
  \State $ \mathbf{z}_{i}^{p} : = \mathbf{z}_{i} + r\cdot \delta \mathbf{z}$; 
  \Statex \textbf{repeat} \Comment{corrector step}
  \State $\mathcal{N}(\mathbf{z}) : = \mathbf{z}_{i}^{p} + \alpha \left[
    \begin{pmatrix}
      \nabla_\mathbf{z} \mathcal{H} (\mathbf{z})\\
      \nabla_\mathbf{z} \mathcal{S} (\mathbf{z},r)^{\top}\\
    \end{pmatrix}^{-1}
    \begin{bmatrix}
      \mathcal{H} (\mathbf{z})\\
      \mathcal{S} (\mathbf{z},r)
    \end{bmatrix}\right]_{\mathbf{z}=\mathbf{z}_i^p}$;
\Statex \textbf{until}~\text{convergence};
\State $    \|  \mathcal{N}^{\kappa}(\mathbf{z}) -    \mathcal{N}^{\kappa-1}(\mathbf{z}) \| \leq \epsilon, \quad \kappa = 1,2,\cdots $
%\Statex $\mathbf{z} := \mathbf{z}_{i+1}$; \Comment{Next point on the curve}
\Statex \textbf{return} $\mathbf{z}_{i+1}$;  \Comment{Next point on the curve}
\State \text{choose a new radius of sphere}~ $r>0$ from \eqref{eq:radius};
\Statex \textbf{until}~\text{traversing is stopped}.
 \EndProcedure
\end{algorithmic}
\end{algorithm}

%%%%%%%%%%%%%%%%%%%%%%%%%%%%%%%%%%%%%%%%%%%%%%%%%%%%%%%%%%%%%%%%%%%%%%%%%%%%%%%%%%%%%%%%%%

\section{Numerical Results}\label{sec:numerics}
This section presents results for two different networks. First, a small three-bus system and later an IEEE 300-bus test case to assess the proposed SC algorithm's computational performance.
\subsection{Three Bus System}
\begin{figure}[!]
    \centering
    \includegraphics[scale=0.07]{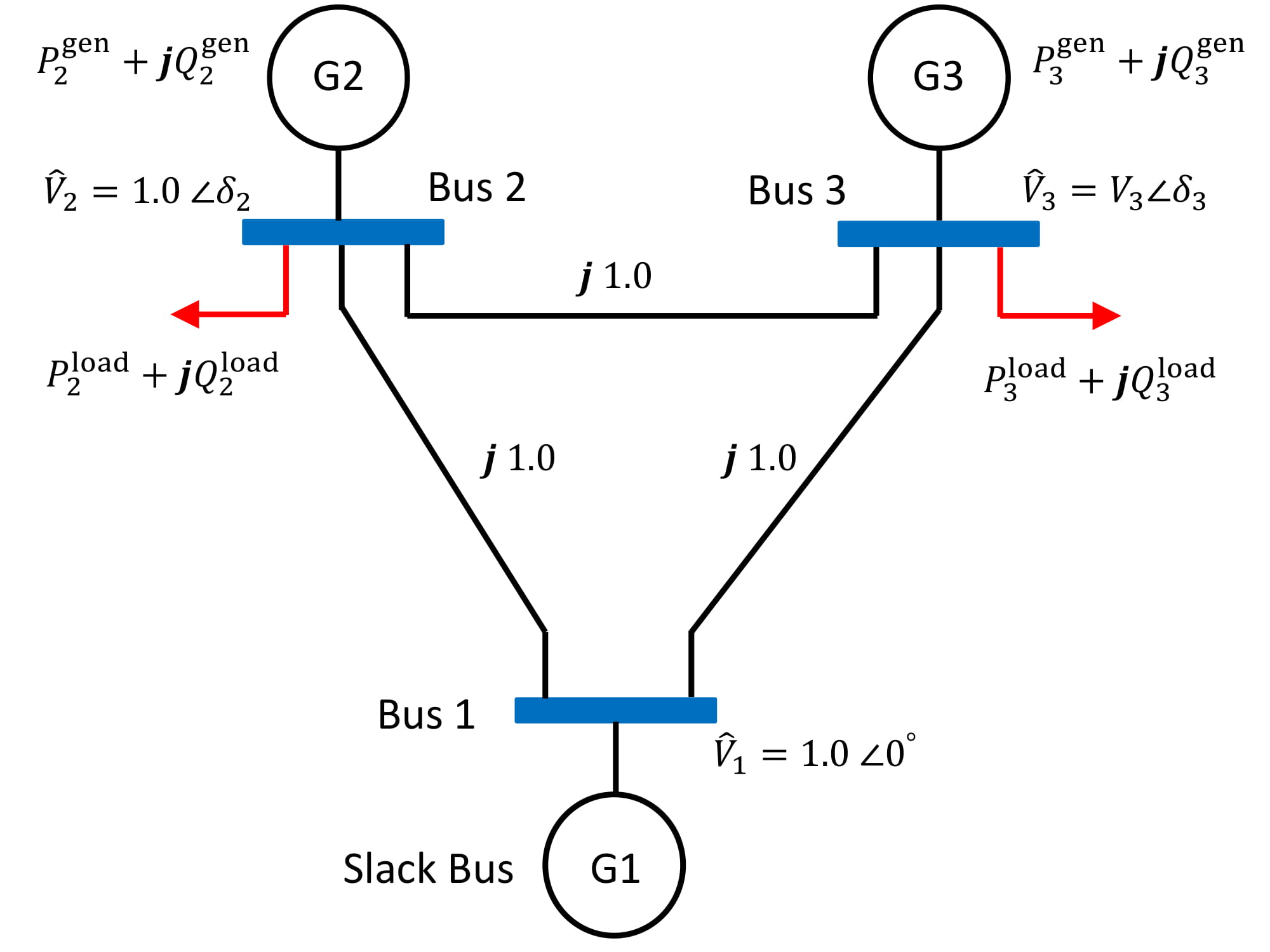}
    \caption{Three bus network.}
    \label{fig:3busd}
\end{figure}

A small three bus system is shown in Figure.\ref{fig:3busd}, which is a slightly modified network from \cite{Hiskens:2001ie}. This network has following settings; the bus 1 is a slack bus with $\hat{V}_{1} = 1.0\angle{0}^{\circ}$. Bus 2 is a PV bus with voltage setpoint $\hat{V}_{2} = 1.0\angle{\delta}_{2}$ and bus 3 is a PQ bus with complex voltage phasor $\hat{V}_{3} = V_{3}\angle{\delta_{3}}$. All the lines in the network are assumed to be lossless with inductance set to $X_{12} = X_{13} = X_{23} = 1.0$ (p.u.). We also assume, that the voltage magnitude at bus $3$ should satisfy \eqref{eq:vlimits}, such that $0.9 \leq |V_{3}| \leq 1.1$ (p.u.).  

\begin{figure*}[t] 
  \centering
  \subcaptionbox{Slice in $P_2 - P_3$ plane.\label{3da1}}[.3\linewidth][c]{%
    \includegraphics[width=.82\linewidth]{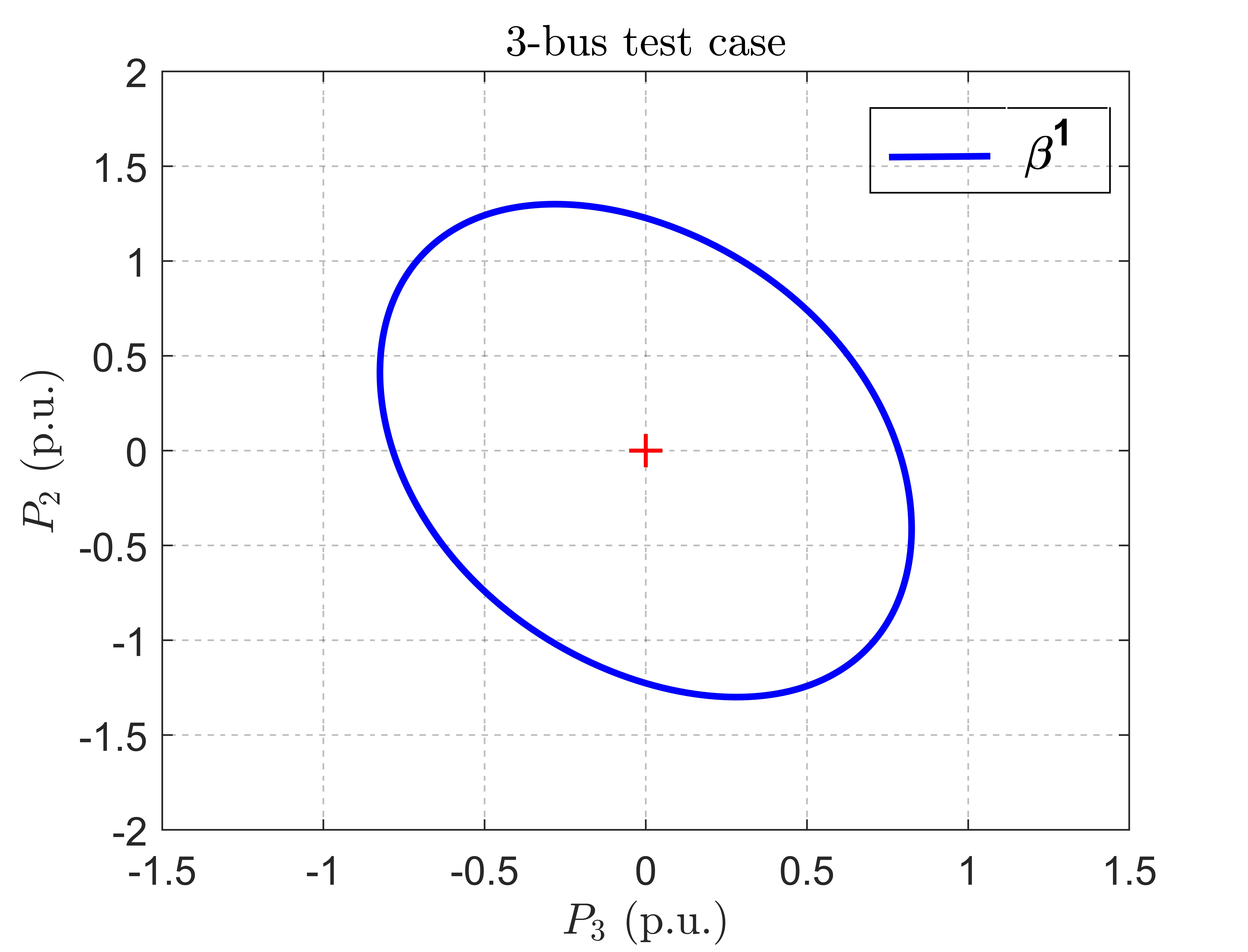}}\quad
  \subcaptionbox{Slice in $P_2 - Q_3$ plane.\label{3da2}}[.3\linewidth][c]{%
    \includegraphics[width=.82\linewidth]{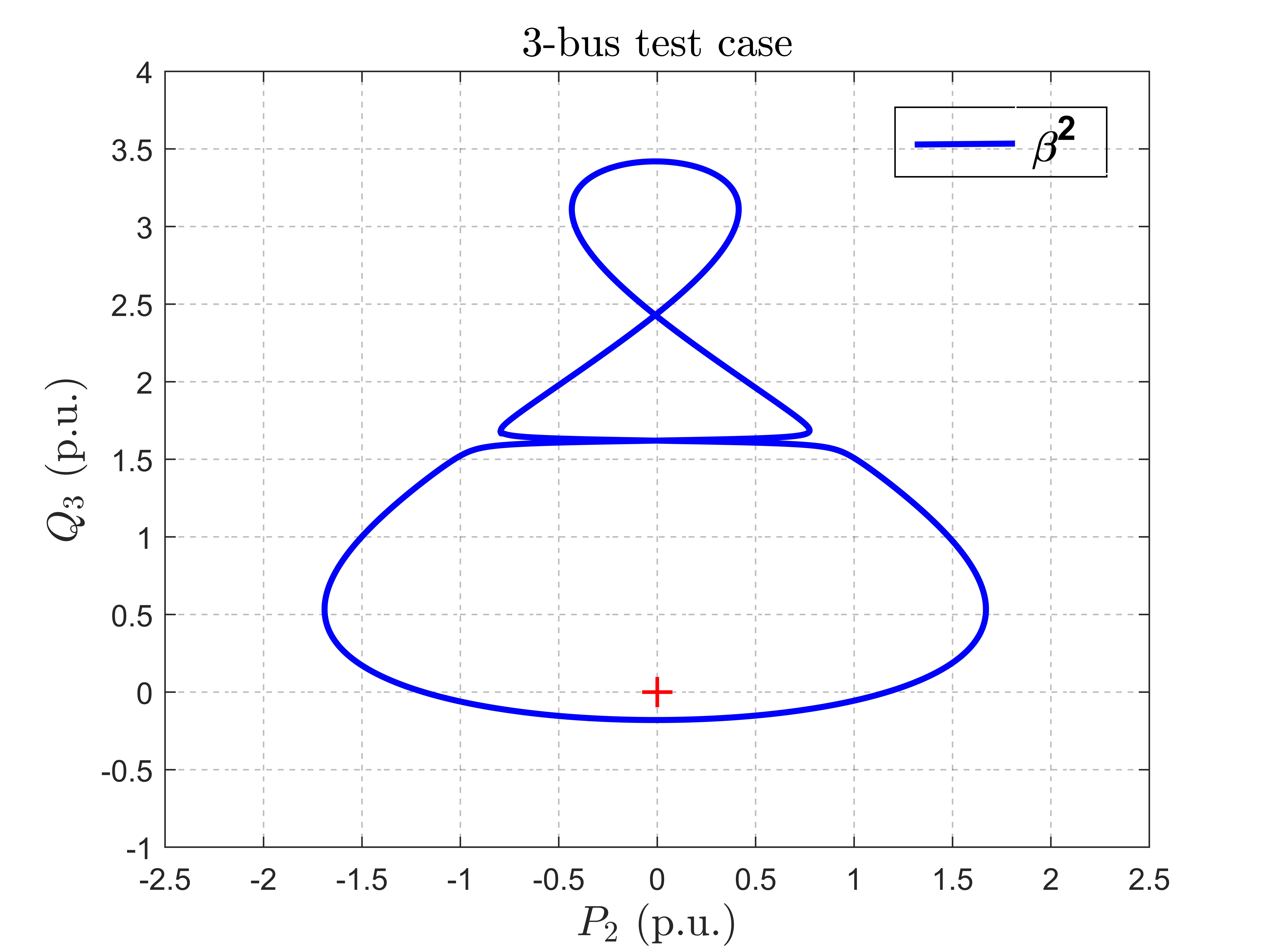}}\quad
  \subcaptionbox{Slice in $P_3 - Q_3$ plane.\label{3da3}}[.3\linewidth][c]{%
    \includegraphics[width=.82\linewidth]{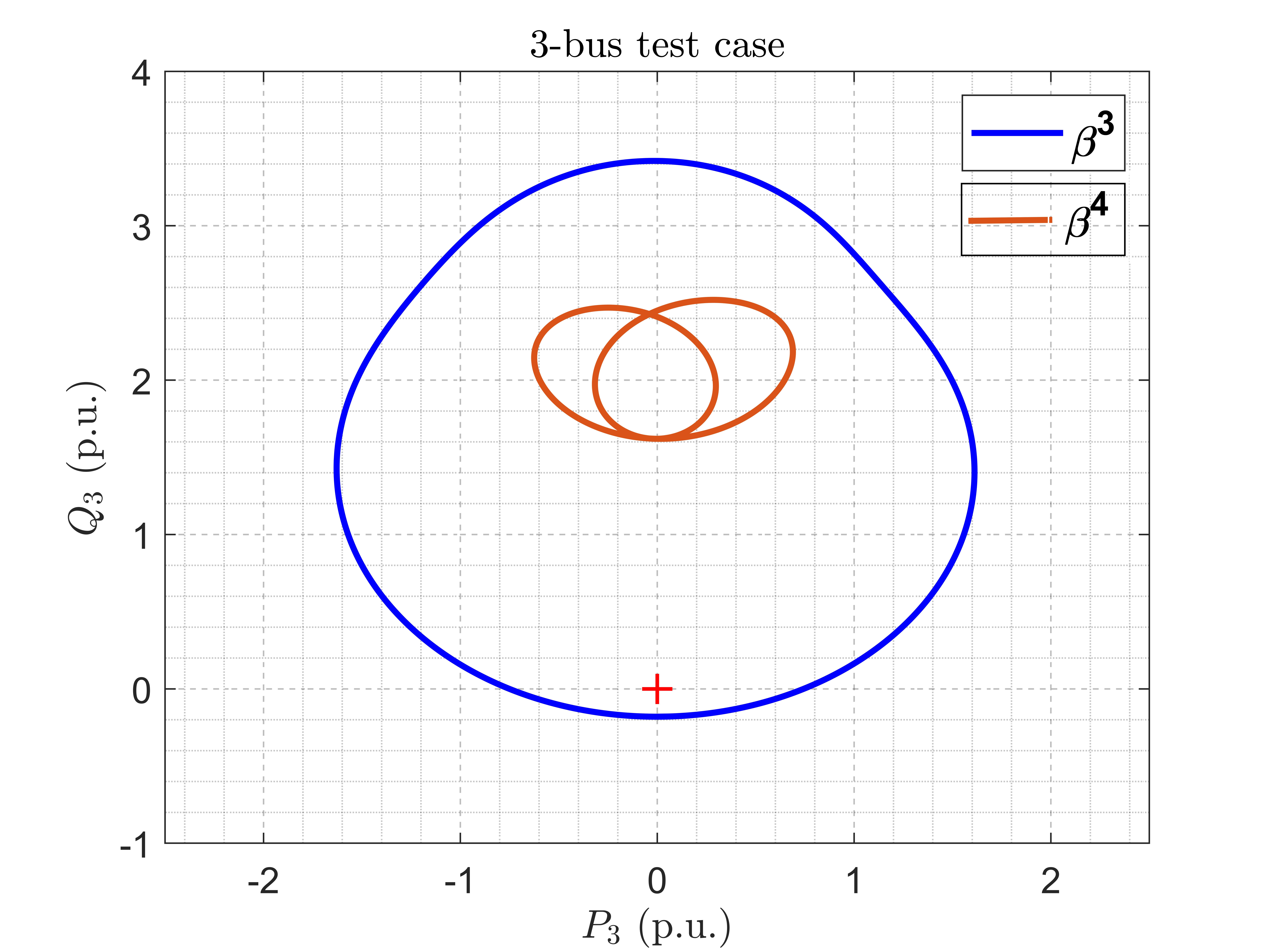}}
   \caption{Cross-sections of the voltage feasibility boundary (3-bus system).\label{fig-3buscross}}
\end{figure*}

\begin{figure*}[h]
  \centering
  \subcaptionbox{Projection in $P_{2}-P_{3}$.\label{3da11}}[.3\linewidth][c]{%
    \includegraphics[width=.82\linewidth]{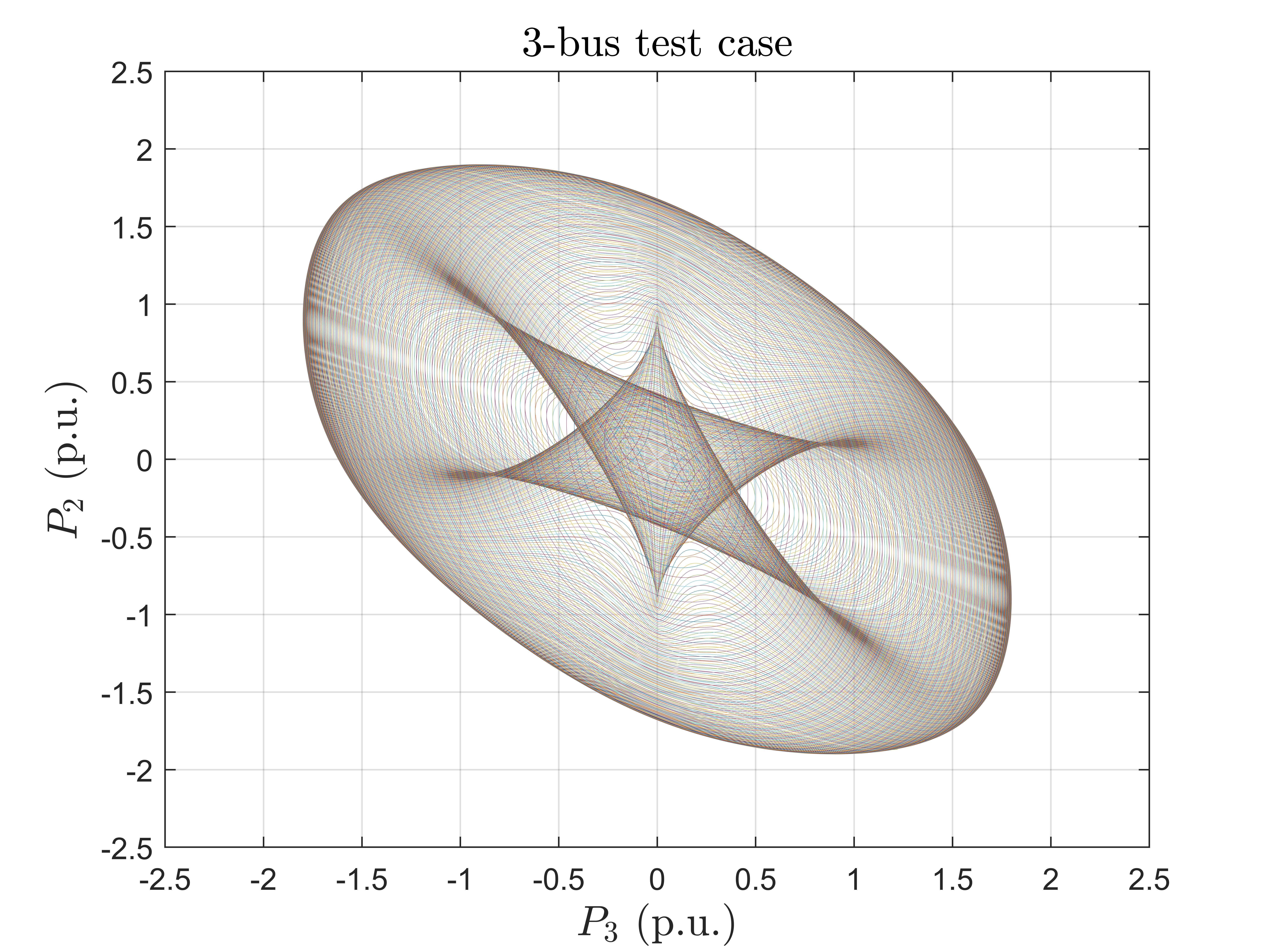}}\quad
  \subcaptionbox{Projection  in $P_{2}-Q_{3}$.\label{3da22}}[.3\linewidth][c]{%
    \includegraphics[width=.82\linewidth]{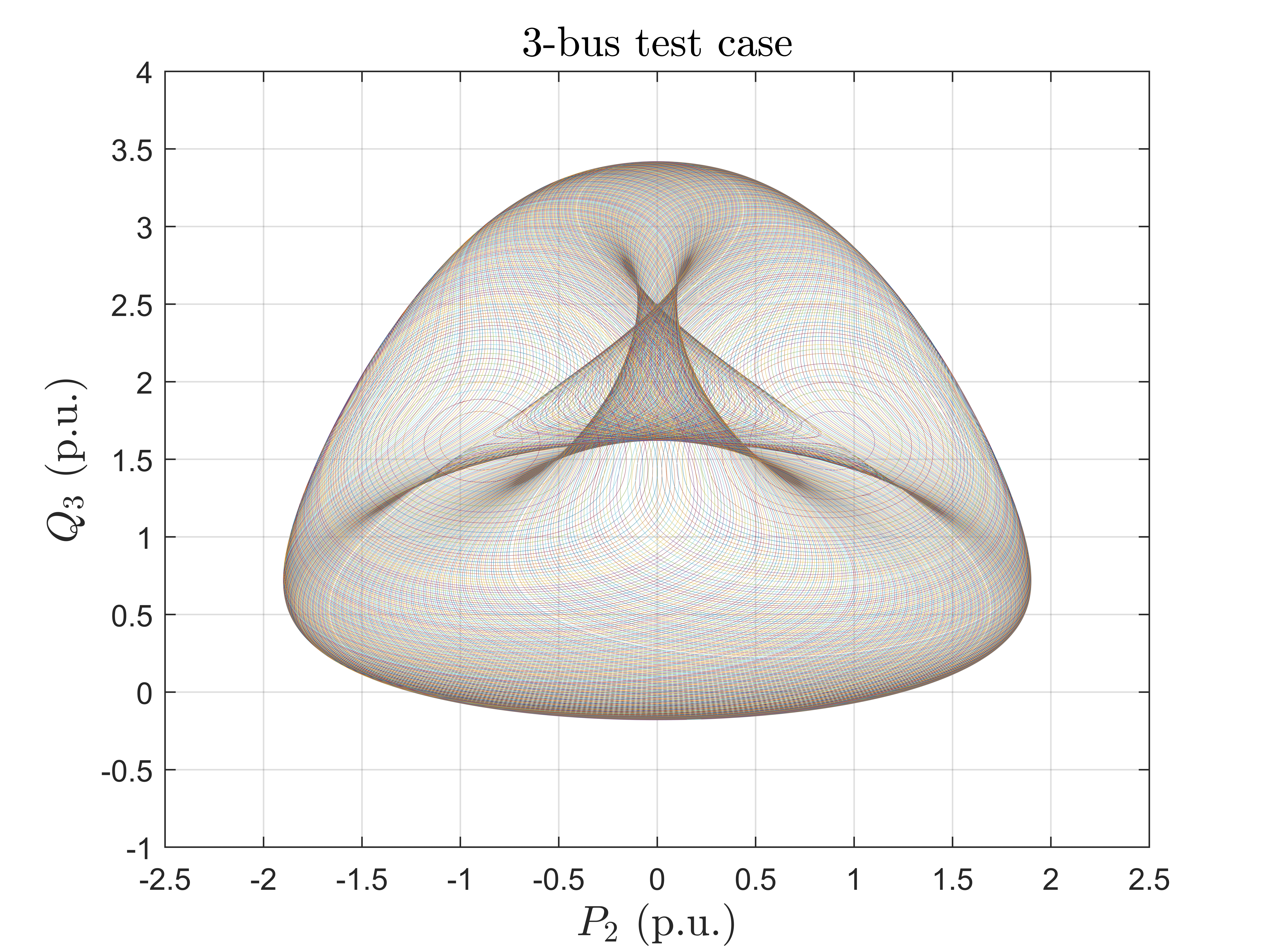}}\quad
  \subcaptionbox{Projection  in $P_{3}-Q_{3}$.\label{3da33}}[.3\linewidth][c]{%
    \includegraphics[width=.82\linewidth]{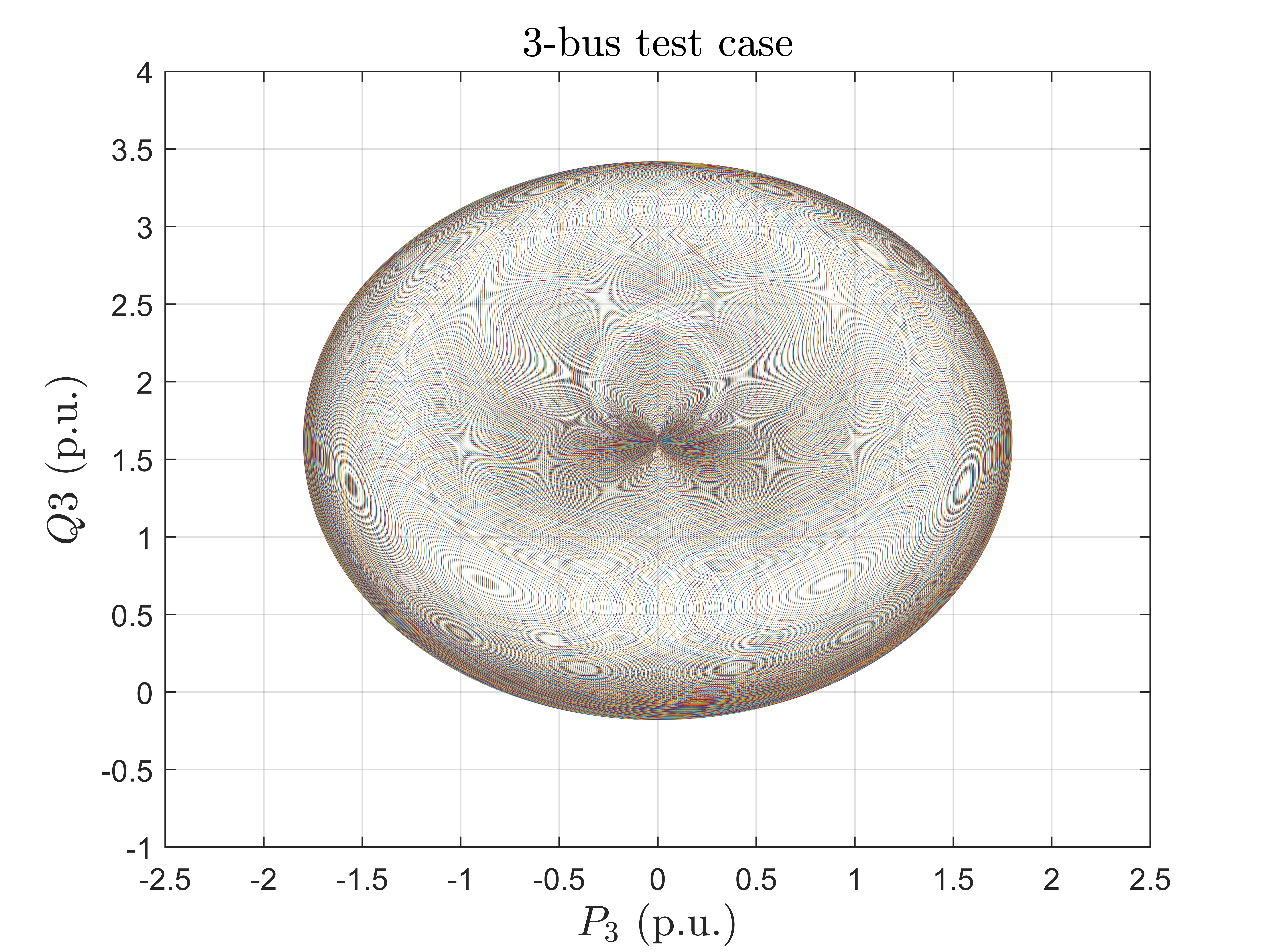}}
   \caption{Projection of the voltage feasibility boundary in $P_{2}-P_{3}-Q_{3}$ (3-bus system).\label{fig-3busproj}}
\end{figure*}

\begin{figure}
    \centering
    \includegraphics[scale=0.015]{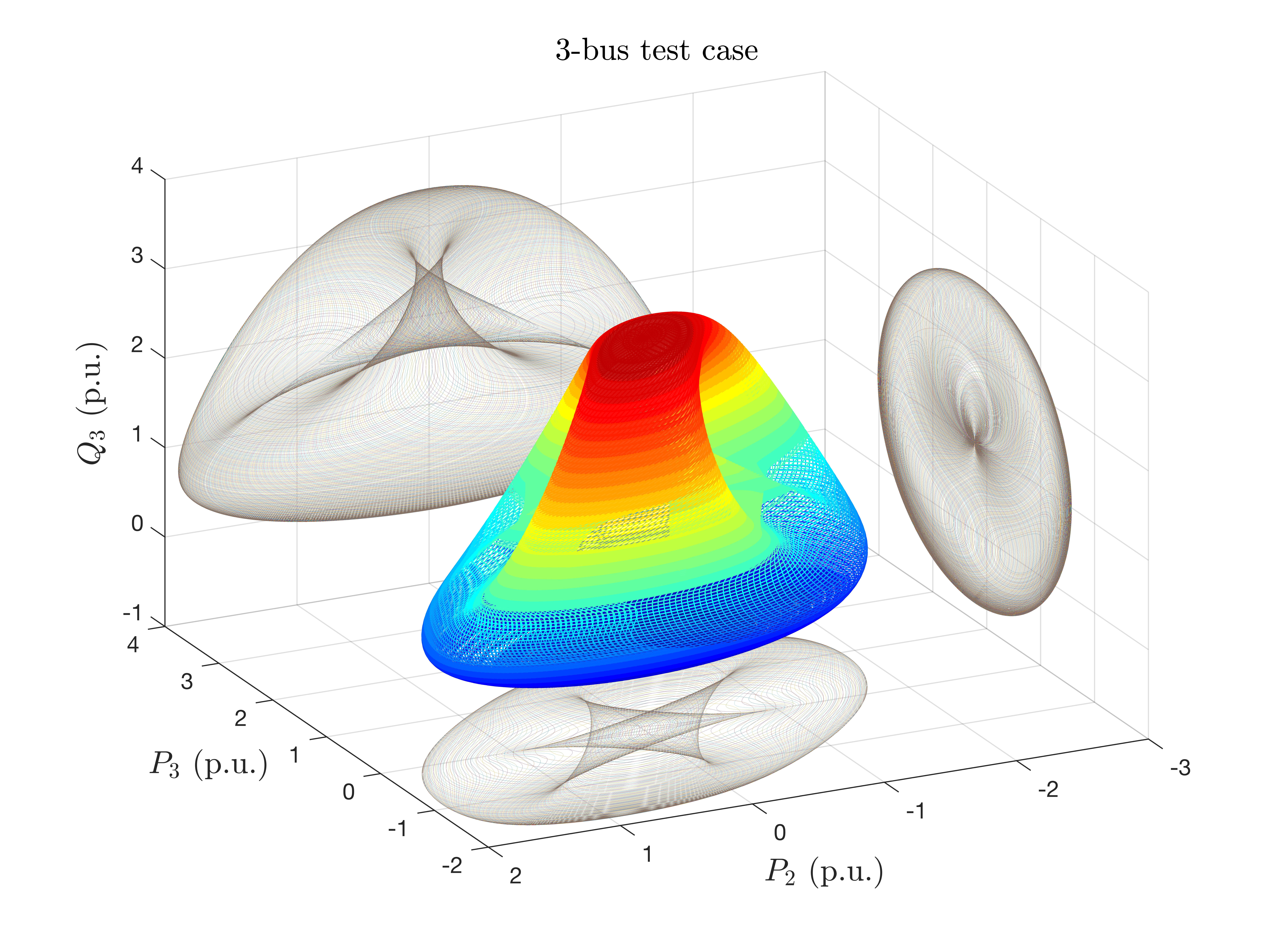}
    \caption{Voltage feasibility boundary surface in $P_{2}-P_{3}-Q_{3}$.} %Also, projection of the feasible manifold is provided in $P_{2}-P_{3}$  as contours of $Q_{3}$. Another projection of the feasible manifold is in $P_{2}-Q_{3}$ plane as contours of $P_{3}$. }
    \label{fig:3xx}
\end{figure}

For the given network topology, the boundary manifold is confined in $P_{2}-P_{3}-Q_{3}$ space. Here, $P_{2}$ and $P_{3}$ describe active power injection at bus 2 and bus 3, respectively. While $Q_{3}$ is the reactive power injection at bus 3. Thus, the feasibility boundary is a certain surface in this $3$-dimensional space. An explicit form of this surface for the given values of parameters is shown in Fig.\ref{fig:3xx}. We note that the shape of the surface is rather complex, even for this very simple three bus network. Next, Fig.\ref{fig-3buscross} illustrates the cross-sections of the surface in three different planes, namely $P_2-P_3$, $P_2-Q_3$, and $P_3 - Q_3$ (every time assuming fixed value of the third parameter), while Fig.\ref{fig-3busproj} shows the projections of the surface to the corresponding planes.  

%Structure of the feasibility space boundary was examined in all possible different cross-sections, namely $P_2-P_3$, $P_2-Q_3$, and $P_3 - Q_3$. Boundary curves in all these cross-sections were computed by applying voltage magnitude limits on bus 3 through \eqref{eq:vlimits} such that $0.9 \leq |V_{3}| \leq 1.1$ (p.u.). %Thus, eliminating the chance to compute boundary curves where inequality on $|V_{3}| $ is violated.
% Also, comparing to proposed method in \cite{Hiskens:2001ie}, the initialization here was rather easy and didn't require an explicit initial guess................... 

%\paragraph{Projection in $P_{2}-P_{3}$} 
%{\color{blue} Unlike in \cite{}, here we formulate the problem with operational constraints reflecting the feasibility problem. Furthermore, the feasibility boundary structure was calculated with $p=3$, resulting in a 2-manifold surface.}

The first cross-section considers the boundary manifold in the $P_2-P_3$ plane. The algorithmic procedure first finds a single-point solution on the $P_2-P_3$ boundary curve for a fixed value of $Q_{3}$. Then the proposed SC algorithm releases the second varying parameter to trace the corresponding 1-dimensional curve. For example, Fig.\ref{3da1} shows a slice of the boundary manifold by $\boldsymbol{\beta}^{1}$ curve in $P_{2}-P_{3}$ plane with $Q_{3} =0.02$ (p.u.). The projection of the manifold in $P_{2}-P_{3}$ space is provided in Fig.\ref{3da11}. Each curve in Fig.\ref{3da11} corresponds to a particular value of $Q_{3}$. From this projection, it is noticeable that the feasibility space boundary is bounded with inner foldings. The space within the inner folds represents the regions where different numbers of equilibria exist. 

%\paragraph{Projection in $P_{2}-Q_{3}$}

Another cross-section examines the feasibility boundary in $P_{2}-Q_{3}$. A similar algorithmic procedure was used to compute  $\boldsymbol{\beta}^{2}$ curve (see Fig.\ref{3da2}) depicting the slice within the $P_{2}-Q_{3}$ plane for $P_3 = 0.02$ (p.u.). It is apparent from Fig.\ref{3da2} that the  $\boldsymbol{\beta}^{2}$ curve is non-convex and also folded. The nomograms of the boundary manifold in $P_{3}-Q_{3}$ space is given in Fig.\ref{3da22} for different values of $P_{3}$. Folding within the boundary manifold is also observable from this projection.    

%\paragraph{Projection in $P_{3}-Q_{3}$}

Finally, the voltage feasibility boundary was considered in $P_{2}-P_{3}$ space. From the previous cases, it was perceived that the solution manifold is clearly folded; hence, there is a possibility to have disconnected boundary curves. Figure.\ref{3da3} depicts a slice in the $P_{3}-Q_{3}$ plane by limiting $P_{2}$ = 0.02 (p.u.). Two disconnected boundary curves, namely $\boldsymbol{\beta}^{3}$ and $\boldsymbol{\beta}^{4}$ are drawn in this slice, as shown in Fig.\ref{3da3}. First, the initial points were calculated corresponding to each of these curves, and then, the proposed SC is initiated to trace these solution boundaries. The nomograms of feasibility boundary in $P_{3}-Q_{3}$ space are given in Fig.\ref{3da33} as contours of $P_{2}$. 

The projections given in $P_2-P_3$, $P_2-Q_3$ and $P_3-Q_3$ were combined to depict the 2-dimensional   boundary surface in $P_{2}-P_{3}-Q_{3}$. Figure. \ref{fig:3xx} shows such surface together with all the projections. In comparison to previous works, this example describes the feasibility space boundary with the operational constraints. Furthermore, the feasibility boundary was calculated with $p=3$, resulting in a 2-dimensional surface. From the computational context, the proposed SC method did not suffer from any numerical stability issues. The variable sphere strategy further improved the computational capabilities, enabling us to compute boundary curves faster and with non-convex sections. Although this three bus network is small, the complexity of the problem is easy to depict. 

% In contrast to the method from \cite{}, the proposed SC does not suffer from initialization issues.  And variable sphere strategy allows us to compute the boundary curves faster with non-convex sections. 

\begin{figure}
    \centering
    \includegraphics[scale=0.013]{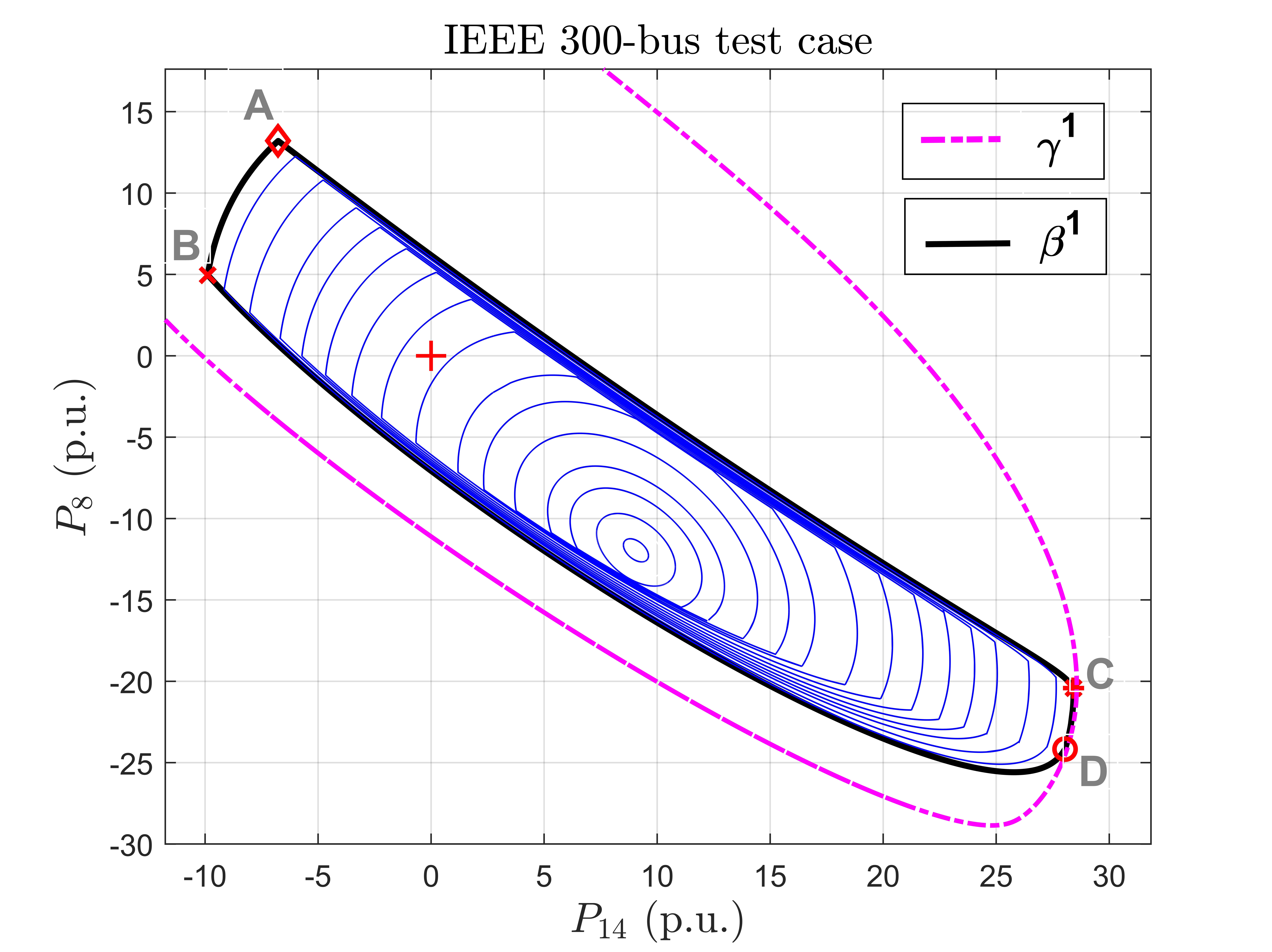}
    \caption{Boundary curves in $P_{8}-P_{14}$ as contours of $Q_{14}$.} %Also, projection of the feasible manifold is provided in $P_{2}-P_{3}$  as contours of $Q_{3}$. Another projection of the feasible manifold is in $P_{2}-Q_{3}$ plane as contours of $P_{3}$.}
    \label{fig:300}
\end{figure}

\subsection{IEEE $300$-bus Test Case}
This section presents results for an IEEE $300$-bus test case to evaluate the proposed SC algorithm's performance in terms of computational tractability and scalability to large networks. The 1-dimensional feasibility boundary curves for this network are shown in Fig.\ref{fig:300}. These curves present the slices of the feasibility boundary in the $P_{8}-P_{14}$ plane for $ 0 \leq Q_{14} \leq 5.8$ (p.u.), and satisfy  the voltage magnitude constraints i.e., $0.9 \leq |V_{i}| \leq 1.1$ (p.u.) on all the $231$ load buses (thus, amounting to $462$ conditions of the type \eqref{eq:vlimits}).

For such a large test network, calculating 1-dimensional boundary curves was a computationally burdensome task. The proposed mathematical formulation with scalar transversality condition allowed the fast computation of these curves without additional computational overhead, even with auxiliary equations representing the inequality constraints. This paper's formulation resulted in an average computational time (for calculating a single curve from Fig.\ref{fig:300}) of about $2.73$ seconds. In contrast, the formulation from \cite{Hiskens:2001ie} resulted in $35.36$ seconds, as it requires solving almost double the number of variables for each point on the boundary.  

The sharp corners are visible in some of the curves from Fig.\ref{fig:300}. For instance, the $\boldsymbol{\beta}^{1}$ curve represents the slice of the feasibility boundary in $P_{8}-P_{14}$ for $Q_{14} = 0$ (p.u.). The corresponding slice of the solvability boundary (i.e., without considering voltage magnitude constraints on load buses) is given by $\boldsymbol{\gamma}^{1}$ curve. Four different sharp turning points namely, ${\mathrm{A,B,C,D}}$ are identified on $\boldsymbol{\beta}^{1}$ curve. These abrupt curvature changes represent the points where different operational constraints are reached. At the section $\mathrm{A-B}, \mathrm{A-C}$ and $\mathrm{B-D}$ the power injections (i.e., $P_8$ and  $P_{145}$) are bounded as the system in \eqref{eq:main} reaches the lower limit of voltage magnitude on bus $14$ and bus $282$. While the section $\mathrm{C-D}$ coincides with the $\boldsymbol{\gamma}^{1}$ manifesting the boundedness of the power injections, i.e., $P_8$ and  $P_{145}$ are due to the solvability condition of the power flow equations. The proposed SC method with adaptive sphere strategy enabled us to trace these curves faster without encountering any numerical instabilities. The traditional PC method from \cite{Hiskens:2001ie} suffered from numerical instabilities around sharp curvature changes as the hyperplane-based continuation becomes slower in the vicinity of sharp turning points, requires a precise adjustment for the continuation step size.

% The proposed SC method with adaptive sphere strategy enabled us to trace these curves faster without encountering any numerical instabilities. The traditional PC method from \cite{Hiskens:2001ie} suffered from numerical instabilities around sharp curvature changes as the hyperplane-based continuation becomes slower in the vicinity of sharp turning points, requires a precise adjustment for the continuation step size.   

For the proposed algorithm, the time required to calculate 1-dimensional boundary curves was compared with the standard method from \cite{Hiskens:2001ie}, and a variant of the SC method from \cite{diaz2017homotopy} to assess the scalability towards large networks. Several runs were performed on different IEEE cases (including the 300-bus test case), and an average time is reported here.  In Table.\ref{table3s3}, parameters $\lambda_{1}$ and $\lambda_{2}$  designates the pair of active power injections for which the 1-dimensional boundary curves were computed. Here $t_{\mathrm{hc}}$, $t_{\mathrm{sc}}$ describe the time for algorithms based on the approach from \cite{Hiskens:2001ie} and \cite{diaz2017homotopy} respectively, while $t_{\mathrm{sc}}^{*}$ delineates the computational time for the proposed SC algorithm with adaptive sphere size strategy. From Table.\ref{table3s3} it is clear that the adaptive SC method introduced in this paper surpasses the traditional approaches from \cite{Hiskens:2001ie} and \cite{diaz2017homotopy}. 

\begin{table}
\caption{Computational time for different IEEE cases.}
\centering
\begin{tabular}{|c|c|c|c|c|c|} 
 \hline
IEEE  &  $\lambda_{1}$ & $\lambda_{2}$ & {$t_\mathrm{hc}$} & {$t_{\mathrm{sc}}$} & {$t_{\mathrm{sc}}^{*}$} \\[0.1 ex]
cases  &   &  & seconds & seconds & seconds\\[0.1 ex]
 \hline\hline
9 Bus     & $P_2$     & $P_{5}$      & 3.21   & 1.66   & 1.10\\[0.5 ex]
14 Bus     & $P_2$     & $P_{3}$     & 5.67   & 1.94   & 1.31\\[0.5 ex]
%30 Bus     & $P_2$     & $P_{3}$    & 5.42   & 2.11   & 1.74\\[0.5 ex]
57 Bus     & $P_8$     & $P_{10}$    & 3.97   & 1.83   & 1.29\\[0.5 ex]
300 Bus    & $P_8$     & $P_{14}$    & 35.36  & 11.1   & 2.73\\[0.5 ex]
 \hline
\end{tabular}
\label{table3s3}
\end{table}

\section{Conclusions and future work}\label{sec:con}

We have introduced an adaptive Spherical Continuation (SC) method with a variable sphere strategy for calculating 1-dimensional boundary curves with computational tractability, speed, and scalability to large networks.  A standardized mathematical formulation is proposed to extend the power flow feasibility problem to a particular equivalent system of the solvability problem, using the slack variable methodology. This approach converts the inequality constraints into a set of polynomial equations suitable to consider in the power flow model, thus preserving the power flow problem's sparsity property.  Here, a scalar transversality condition was used to enforce the solution of the power flow equations on the boundary of feasibility space, which allowed computing the boundary points without initialization heuristics, smaller computational overhead, and better scalability.

We show that even for a small 3-bus system, the feasibility boundary is a complex manifold with folded parts representing multiple equilibrium regions. The SC algorithm's performance was assessed on an IEEE 300-bus network, calculating 1-dimensional boundary curves faster than the existing algorithms without encountering any numerical instabilities. We also show that the SC algorithm can trace section of the curves with sharp corners or non-convexity. Future work includes using the method for potential practical applications like contingency analysis, improving local optima in relaxed OPF problems or determining an electric grid's transfer capabilities, etc.

\ifCLASSOPTIONcaptionsoff
  \newpage
\fi

\bibliographystyle{IEEEtran}
\bibliography{biblography.bib}

% that's all folks
\end{document}